%% file: IRREG.TEX
\documentclass[12pt]{report}
\usepackage{amsfonts, amsmath, amssymb}
\newcommand{\rr}{\mathbb R}

\newcommand{\zz}{\mathbb Z}
\newcommand{\cc}{\mathbb C}
\renewcommand{\ggg}{\mathbb G}

\newcommand{\RR}{{\mathfrak R}}

\newcommand{\FF}{{\mathfrak F}}

\newcommand{\LL}{\mathfrak L}

\newcommand{\CC}{{\mathfrak C}}
\newcommand{\SC}{{\mathfrak S}{\mathfrak C}}

\newtheorem{theorem}{Theorem}
\newtheorem{lemma}{Lemma}
\newtheorem{prop}{Proposition}
\newtheorem{rem}{Remark}
\newtheorem{cor}{Corollary}

\newtheorem{exmp}{Example}

\title{ Irregular subsets of the Grassmannian manifolds}
\author{M. A. Pankov}
\date{}

\begin{document}

\maketitle

\tableofcontents

\input{INTRODUC.TEX}
\input{CHAPTER1.TEX}
\input{CHAPTER2.TEX}
\input{CHAPTER3.TEX}
\input{BIB.TEX}

\end{document}

%% file: INTRODUC.TEX
\chapter*{Introduction}

\label{intro}
\addtocontents{toc}{\contentsline{chapter}{Introduction}{\pageref{intro}}}

Irregular subsets of the Grassmannian manifolds were
introduced in the author's papers [P1, P2]
where some results connected with the Chogoshvili Conjecture
were optained.
The Conjecture states that for any $k$-dimensional
compact set $X\subset \rr^{n}$
(a set is called $k$-dimensional if its topological
dimension is equal to $k$ [Eng])
there exists a
$k$-dimensional plane such that the orthogonal projection
of $X$ onto this plane is a stable mapping.
We do not define this notion here.
The exact formulation of the Chogoshvili Conjecture
could be found in [Dr] (see also the original
Chogoshvili's paper [Chog]).
It is well-know that the Chogoshvili
Conjecture fails for non-compact $k$-dimensional sets
(the respective example was constructed by
K. Sitnikov [S, Eng]). For the compact case
the Conjecture was opened for a long time.
Recently A. N. Dranishnikov [Dr] disproved
it.

In the papers [P1, P2] for arbitrary set $X\subset \rr^{n}$
we considered the set $R^{n}_{k}(X)$ of $k$-dimensional
planes in $\rr^{n}$ such that the orthogonal projections
of $X$ onto these planes are "regular"
(note that if some projection is a stable mapping
then it is regular, the inverse statement does not hold).
It was proved there that if our set $X$ satisfies
the condition $\dim X \ge k$ then
$$\ggg^{n}_{k}\setminus R^{n}_{k}(X)$$
is an irregular subset of the Grassmannian
manifold $\ggg^{n}_{k}$.
It seems to be natural to ask
how large may be an irregular set?
This problem was studied in [P2].
If $k=1,n-1$ then
each irregular subset of
$\ggg^{n}_{k}$ is nowhere dense.
For the general case
($1<k<n-1$) the similar statement is not proved.  However, there
are a few results supporting our conjecture.  One of them states
that any irregular subset  has an empty interior
in $\ggg^{n}_{k}$.

This book is devoted to study
general properties of irregular sets and is not connected
with Dimension Theory.

Let $V$ be an $n$-dimensional vector space under
arbitrary field and
$\ggg^{n}_{k}$ be
the Grassmannian manifold of $k$-dimensional
linear subspaces of $V$.
The definition of irregular sets are based on the
notion of a regular subset of the Grassmannian
manifold.

A subset of $\ggg^{n}_{k}$ is called {\em
regular} if
planes belonging to it
are coordinate planes for some coordinate system
for $V$.
A regular subset of $\ggg^{n}_{k}$ is a discrete set
containing at
most $c^{n}_{k}$ elements,
here $c^{n}_{k}$ is the respective binomial
coefficient (a coordinate system for $V$
has $c^{n}_{k}$ distinct
$k$-dimensional coordinate planes).
A regular subset of $\ggg^{n}_{k}$
containing $c^{n}_{k}$ elements will be called
{ \it maximal}.
Regular sets could be consider as some generalization
of collections of linearly independent lines.

We say that a subset of $\ggg^{n}_{k}$ is
{\it irregular} if it is not regular and does
not contain maximal regular subsets.
Each irregular set is contained in
some maximal irregular set
(an irregular set  is maximal if any
irregular set containing it coincides with it).

There are a few simple examples of irregular subsets of
$\ggg^{n}_{k}$.

\begin{enumerate}
\item[---] The cases $k=1,n-1$ are trivial.
Each maximal irregular subset
of $\ggg^{n}_{1}$ consists of all lines
contained in some $(n-1)$-dimensional plane.
Each maximal irregular subset
of $\ggg^{n}_{n-1}$ consists of all
planes containing a fixed line.
\item[---] The general case is more complicated.
Let $s\in \ggg^{n}_{m}$ and $m\ne k$.
Consider the set of all
planes $l\in \ggg^{n}_{k}$ satisfying the
condition
$$\dim s\cap l \ge 1\;.$$
In the case $m\le n-k$ our set is irregular
(if $m>n-k$ then it coincides with $\ggg^{n}_{k}$).
Moreover, in the case $m=n-k$
this irregular set is maximal.
\end{enumerate}

A transformation $f$ of $\ggg^{n}_{k}$
(bijection onto itself)
is called {\it regular} if $f$ and $f^{-1}$
map any regular set onto regular set.
It is trivial that $f$ is regular
if and only if $f$ and $f^{-1}$
preserve the class of irregular sets.

Consider two examples of regular transformations.
\begin{enumerate}
\item[---] Each linear transformation
of $V$ defines a regular transformation of $\ggg^{n}_{k}$.
\item[---] Each non-singular bilinear form $\Omega$ on $V$
defines the bijection
of $\ggg^{n}_{k}$  onto $\ggg^{n}_{n-k}$
transferring  each plane $s\in \ggg^{n}_{k}$
to the $\Omega$-orthogonal complement.
This bijection maps any regular subset of $\ggg^{n}_{k}$ onto
a regular subset of $\ggg^{n}_{n-k}$.
For the case $n=2k$ it is a regular transformation
of $\ggg^{2k}_{k}$.
\end{enumerate}

One of the main results of the book states that
\begin{enumerate}
\item[---] if $n\ge 3$ and $n\ne 2k$ then
regular transformations of $\ggg^{n}_{k}$
are defined by linear transformations;
\item[---] if $k\ge 2$ and $n=2k$ then
regular transformations of $\ggg^{2k}_{k}$
are defined by linear transformations
or non-singular bilinear forms.
\end{enumerate}
In the case $k=1$ this statement is known as
the Fundamental Theorem of Projective Geometry [Art, O'M, Die].
For the general case it was proved by author
[P3].

We say that two subsets
of $\ggg^{n}_{k}$
are {\it similar} if there exists
a regular transformation of $\ggg^{n}_{k}$
transferring one of this sets to other set.
In the cases $k=1,n-1$ any two maximal irregular subsets
of $\ggg^{n}_{k}$ are similar. For the general case
it fails.
It was proved above that
for any two natural numbers
$n$ and $k$ satisfying the conditions
$n>3$ and $1<k<n-1$ there exist
maximal irregular subsets
of $\ggg^{n}_{k}$ which are not
similar.

It seems to be natural to ask how many non-similar
maximal irregular sets exist? This problem is very
difficult and we can not consider it here.

{\bf Thanks:}
I wish to express my deep gratitude to my teacher
V. V. Sharko and to my colleagues
I. Yu. Vlasenko, S. I. Maksimenko, E. A. Polulyakh
for their interest in our research.
I also want to thank my wife and mother  who have boundlessly
supported me while I have been writing this book.

%% file: CHAPTER1.TEX
\chapter{Grassmannian manifolds and their transformations}

First of all (Section 1.1) we recall some standart facts from
Linear Algebra connecting with linear mappings and bilinear
forms.

Denote by $\ggg^{n}_{k}(V)$
the Grassmannian manifold of $k$-dimensional
linear subspaces of an $n$-dimensional vector space
$V$.
In Section 1.2
the following objects will be considered:
\begin{enumerate}
\item[---] transformations of $\ggg^{n}_{k}(V)$
defined by linear
transformations of $V$;
\item[---] bijections of $\ggg^{n}_{k}(V)$
onto $\ggg^{n}_{n-k}(V)$ defined by
non-singular bilinear forms on $V$;
\item[---] transformations of $\ggg^{n}_{m}(V)$
induced by
transformations of $\ggg^{n}_{k}(V)$.
\end{enumerate}
The term "transformation" means a bijection
onto itself.

Sections 1.3 and 1.4 are devoted to prove the Fundamental Theorem
of Projective Geometry and the Chow Theorem on
transformations of the Grassmannian manifolds
preserving the distance between
planes.

All main results of this chapter, except Theorem 1.2.1,
were taken from [Art, O'M, Die].

\section{Linear mappings and bilinear forms}
\setcounter{equation}{0}
\setcounter{prop}{0}
\setcounter{lemma}{0}
\setcounter{defn}{0}

\subsection{Linear transformations of vector spaces}

Let $F$ be a field. The composition
operation defines on the class of automorphisms of $F$
the natural group structure.
In what follows this group will be denoted
by $Aut(F)$.

Let us consider a few examples.

\begin{exmp}{\rm
$Aut(\rr)=\{Id\}$.
$\blacksquare$
}\end{exmp}

\begin{exmp}{\rm
$Aut(\cc)=\{Id, J\}$, where
$J(z)=\overline{z}$ for each $z\in \cc$
($\overline{z}$ is the complex conjugate to $z$).
$\blacksquare$
}\end{exmp}

Let $V$ and $V'$ be  finite-dimensional vector spaces under
the field $F$ and $\sigma \in Aut(F)$.
We say
that a mapping $f: V\to V'$ is $\sigma$-{\it linear} if
it satisfies the next conditions:
$$f(x+y)=f(x)+f(y)\;\;\;\forall\;x,y\in V$$
and
$$ f(ax)=\sigma(a)f(x)\;\;\;\forall\;x\in V,a \in F\;.$$
In the case
when
$$\dim V=\dim V' \mbox{ and } Ker(f)=\{0\}$$
the mapping $f$ is called a $\sigma$-{\it linear
isomorphism} of $V$ onto $V'$.

Denote by $\LL_{\sigma}(V,V')$ the
class of $\sigma$-linear isomorphisms of $V$ onto $V'$ and
define
$$\LL(V,V')=\bigcup_{\sigma \in Aut(F)}\LL_{\sigma}(V,V')\;.$$
In the case when $V=V'$
we write $\LL(V)$ and  $\LL_{\sigma}(V)$
in place of
$\LL(V,V')$ and $\LL_{\sigma}(V,V')$,
respectively.
Elements of
$\LL(V)$ will be called {\it linear transformations} of $V$.

\begin{exmp}{\rm
Let $\sigma \in Aut(F)$ and $B=\{x_{i}\}^{n}_{i=1}$
be a base for $V$.
Consider the linear transformation
$f_{\sigma B}\in \LL_{\sigma}(V)$ defined by the condition
$$f_{\sigma B}(x_{i})=x_{i}\;\;\;\forall\;i=1,...,n\;.$$
For any $f\in \LL_{\sigma}(V)$ there exists
unique $g\in \LL_{Id}(V)$ such that
$$f(x_{i})=g(x_{i})\;\;\;\forall\;i=1,...,n\;.$$
The linear transformation $f$ could be represented as the
composition $f=gf_{\sigma B}$.
$\blacksquare$ }\end{exmp}

\begin{prop}
The composition operation
defines on
$\LL(V)$  the group structure such that
$\LL_{Id}(V)$ is a normal subgroup of $\LL(V)$,
$$\LL(V)/ \LL_{Id}(V) \cong Aut(F)$$
and each factor class coincides with some
$\LL_{\sigma}(V)$.
\end{prop}

{\bf Proof.}
First statement is a trivial consequence of the
following equations
\begin{equation}
f_{2}f_{1}\in \LL_{\sigma_{2}\sigma_{1}}(V)\;\;\;
\forall\;f_{1}\in \LL_{\sigma_{1}}(V),
f_{2}\in \LL_{\sigma_{2}}(V)\;,
\end{equation}
\begin{equation}
f^{-1}\in \LL_{\sigma^{-1}}(V)\;\;\;
\forall\;f\in \LL_{\sigma}(V)
\end{equation}
showing also  that
$$f\LL_{Id}(V)f^{-1}\subset\LL_{Id}(V)\;\;\;
\forall\;f\in \LL(V)\;.$$
Therefore, the subgroup $\LL_{Id}(V)$ is normal.

Two linear transformations $f$ and $g$ of $V$
belong to  same factor class if and
only if $f^{-1}g \in \LL_{Id}(V)$.
Equation (1.1.1) and (1.1.2) show that
the last condition is equivalent
to the  existence of $\sigma \in Aut(F)$
such that $f,g \in \LL_{\sigma}(V)$.
This implies that
the factor class  containing $f$ and $g$
coincides with $\LL_{\sigma}(V)$.

Equations (1.1.1) and (1.1.2)
guarantee
that the one-to-one correspondence
$$\LL_{\sigma}(V)\longleftrightarrow \sigma$$
is an isomorphism between the groups
$\LL(V)/ \LL_{Id}(V)$ and $Aut(F)$.
$\blacksquare$

\begin{rem}{\rm
Each linear isomorphism $f\in \LL_{\sigma}(V,V')$
defines the isomorphism
$$g\to  fgf^{-1}$$
of the group $\LL(V)$
onto the group $\LL(V')$ transferring $\LL_{Id}(V)$
to $\LL_{Id}(V')$.
$\blacksquare$
}\end{rem}

\subsection{Dual space}

Let  $\sigma \in Aut(F)$.
Denote by $V^{*}_{\sigma}$
the class of $\sigma$-linear functionals on $V$.
It has the natural structure of a vector space
and $\dim V^{*}_{\sigma}=\dim V$.
In the case when $\sigma =Id$  this class will be denoted by
$V^{*}$
and will be called
the {\it dual} space.

\begin{rem}{\rm
Each automorphism $\sigma \in Aut(F)$ defines
the $\sigma$-linear isomorphism
$$f\to \sigma f$$
of $V^{*}$ onto $V^{*}_{\sigma}$.
$\blacksquare$
}\end{rem}

Let $B=\{x_{i}\}^{n}_{i=1}$
be a base for $V$.
For each $i=1,...,n$
there exists unique $x^{i} \in V^{*}$
such  that
$$x^{i}(x_{i})=1 \mbox{ and } x^{i}(x_{j})=0
\mbox{ if } i\ne j\;.$$
Then $B^{*}=\{x^{i}\}^{n}_{i=1}$
is a base for the space $V^{*}$.
We say that this base is
{\it dual} to $B$.
The mapping
$$\Sigma^{n}_{i=1} a_{i}x_{i}\to
\Sigma^{n}_{i=1}a_{i}
x^{i}$$
is a linear isomorphism of $V$
onto the dual space $V^{*}$.
This isomorphism is not canonical,
because it depends on $B$.

Any $x \in V$ could be considered as the linear functional
on $V^{*}$ transferring each $f\in V^{*}$ to $f(x)$; i.e.
$V \subset V^{**}$.  Inversely, each linear functional on $V^{*}$
could be represented in this form
and the imbedding $V \subset V^{**}$ is a linear isomorphism
of $V$ onto $V^{**}$.
This isomorphism is canonical
and we can write $V =V^{**}$.

For a linear subspace $U \subset V$
the set
$$U^{\perp}=\{\;f\in V^{*}\;|\;f(x)=0\;\;\forall\;x\in U\;\}$$
is a linear subspace of
$V^{*}$ the dimension of which is equal to
$\dim V -\dim U$.
This subspace is called the {\it orthogonal complement} to
$U$.
It is not difficult to see that $U^{\perp \perp}=U$.

\subsection{Bilinear forms}

Let $\sigma_{1},\sigma_{2}\in Aut(F)$.
We say that a mapping
$$\Omega :V\times V \to F$$
is a
$(\sigma_{1}, \sigma_{2})$-{\it bilinear form} on $V$
if
$$\Omega_{1}(x)=\Omega(x,\cdot)\in V^{*}_{\sigma_{2}}$$
and
$$\Omega_{2}(x)=\Omega(\cdot ,x) \in V^{*}_{\sigma_{1}}$$
for each $x\in V$.
In other words, $\Omega_{1}$ and $\Omega_{2}$
are $\sigma_{1}$-linear and $\sigma_{2}$-linear
mappings of $V$ into the spaces
$V^{*}_{\sigma_{2}}$ and $V^{*}_{\sigma_{1}}$,
respectively.

\begin{exmp}{\rm
Let $\{x_{i}\}^{n}_{i=1}$ be a base for
$V$ and $\sigma_{1},\sigma_{2}\in Aut(F)$.
For any
two vectors
$x=\Sigma^{n}_{i=1} a_{i}x_{i}$
and
$y=\Sigma^{n}_{i=1} b_{i}x_{i}$
define
$$\Omega(x,y)=
\Sigma^{n}_{i=1} \sigma_{1}(a_{i})\sigma_{2}(b_{i})\;.$$
Then $\Omega$ is a $(\sigma_{1},\sigma_{2})$-bilinear form.
$\blacksquare$
}\end{exmp}

The  $(\sigma_{2}, \sigma_{1})$-bilinear
form $\Omega'$  defined by the condition
$$\Omega'(x,y)=\Omega(y,x)\;\;\;\forall\;x,y\in V$$
will be called {\it inverse}
({\it conjugate})
to $\Omega$.
It is trivial that
$\Omega'_{1}=\Omega_{2}$
and $\Omega'_{2}=\Omega_{1}$.

Now fix some  base $B=\{x_{i}\}^{n}_{i=1}$ for $V$
and consider the matrix
$$M=\{\Omega(x_{i},x_{j})\}_{i\,j=}%
{\vphantom{f({\bf x}_{i},{\bf x}_{j})}}_{1}^{n}\;.$$
Then
$$\Omega(x,y)=\Sigma_{i\,j=}
{\vphantom{\Sigma}}_{1}^{n}\Omega(x_{i},x_{j})
\sigma_{1}(a_{i}) \sigma_{2}(b_{j})
\;\;\;\forall\;
x=\Sigma^{n}_{i=1} a_{i}x_{i},
y=\Sigma^{n}_{i=1} b_{i}x_{i}
\;.$$
This implies that the matrixes of
the linear mappings
$\Omega_{1}$ and $\Omega_{2}$
in the base $B$
coincide with $M$ and $M^{t}$, respectively
($M^{t}$ is the matrix transponed to $M$).

Our $(\sigma_{1}, \sigma_{2})$-bilinear
form $\Omega$ is called
{\it non-singular} if
the mapping $\Omega_{1}$ is an isomorphism.
It holds if and only if
$\Omega_{2}$ is an isomorphism.
Therefore, $\Omega$ is non-singular
if and only if the inverse form
is non-singular.

{\bf Each bilinear form considering in what follows will be
non-sigular!!!}

A linear transformation $f \in\LL(V)$
defines a bijection $f^{*}$ of the class
of bilinear forms on $V$ onto
itself.
For each bilinear form $\Omega$  the form
$f^{*}(\Omega)$ is defined by the following condition
$$f^{*}(\Omega)(x,y)=\Omega(f(x),
f(y))\;\;\;\forall\;x,y\in V\;.$$
Two bilinear forms
$\Omega'$ and $\Omega''$ on $V$
are {\it equivalent} (or {\it similar})
if there exists
a linear transformation $f \in\LL(V)$
such that $\Omega''=f^{*}(\Omega')$.

\subsection{Reflexive forms}

We say that some bilinear form $\Omega$ on $V$
is  {\it reflexive} if for any two vectors $x,y\in V$ the
conditions
$$\Omega(x,y)=0 \mbox{ and }\Omega(y,x)=0$$
are
equivalent.

Let $\Omega$ be a bilinear form on $V$.
For each linear subspace $U \subset V$
the set
$$U^{\perp}_{\Omega}=\{\;y\in V\;|\;\Omega(x,y)=0\;\;\forall\;x\in
U\;\}$$
is a linear subspace of $V$ the dimension of which
is equal to $\dim V-\dim U$. It is called the
$\Omega$-{\it orthogonal complement} to $U$.
Our form $\Omega$ is reflexive if and only if
the equality
$$(U^{\perp}_{\Omega})^{\perp}_{\Omega}=U$$
holds for any linear
subspace $U \subset V$.

Consider a few examples of reflexive forms.

\begin{exmp}{\rm
The form $\Omega$ considered in Example
1.1.4 is reflexive if
$\sigma_{1}=\sigma_{2}$.
Therefore,
on each finite-dimensional vector space  there exists
a reflexive bilinear form.
}
$\blacksquare$
\end{exmp}

\begin{exmp}{\rm
Let $\Omega$ be a bilinear form on $V$
satisfying the condition
$$\Omega(x,y)=r\Omega(y,x)\;\;\;\forall\; x,y\in V\;.$$
Then it is reflexive.
If $r=1$ then $\Omega$
is called {\it symmetric},
in the case $r=-1$ it will be called {\it skew-symmetric}.
$\blacksquare$
}\end{exmp}

\begin{exmp}{\rm
We say that an $(Id,\sigma)$-bilinear
form $\Omega$ on $V$ ($\sigma \ne Id$)
is {\it hermitian} if
$$\Omega(x,y)=\sigma(\Omega(y,x))\;\;\;\forall\;x,y\in V$$
and $\sigma$ is an involution (i.e. $\sigma^{2} =Id$).
Each hermitian form is reflexive.
Moreover, for any $a\in F\setminus \{0\}$
the form $\Omega'=a\Omega$ is reflexive,
but it is not hermitian.
Consider, for example, $a\in F$
such that $\sigma(a)=-a$.
In this case we have
$$\Omega'(x,y)=-\sigma(\Omega'(y,x))\;\;\;\forall\;x,y\in V\;.$$
An $(Id,\sigma)$-bilinear
form $\Omega'$ on $V$ ($\sigma$ is an involution)
satisfying the
last condition is called {\it skew-hermitian}.
$\blacksquare$
}\end{exmp}

\begin{exmp}{\rm
Let $F=\cc$ and $\{x_{i}\}^{n}_{i=1}$
be a base for $\cc^{n}$. Define
$$\Omega(x,y)=
\Sigma^{n}_{i=1} a_{i}\overline{b_{i}}
\;\;\;
\forall\;x=\Sigma^{n}_{i=1} a_{i}x_{i},
y=\Sigma^{n}_{i=1} b_{i}x_{i}\;.$$
Then $\Omega$ and $i\Omega$ are hermitian and skew-hermitian
$(Id, J)$-bilinear forms, respectively.
$\blacksquare$
}\end{exmp}

\begin{rem}{\rm
Let $\Omega$
be a $(\sigma_{1},\sigma_{2})$-bilinear
form. Then $\sigma^{-1}_{1}\Omega$ is an
$(Id,\sigma^{-1}_{1}\sigma_{2})$\\-bilinear
form. It
is reflexive if and only if the form $\Omega$ is
reflexive. Reflexive
$(Id,\sigma)$-bilinear forms
were  classified by
G. Birkhoff and J. von Heumann [BH].
Their result states that
each reflexive $(Id,\sigma)$-bilinear
form $\Omega$ satisfies one of the following conditions:
\begin{enumerate}
\item[---] $\sigma= Id$ and $\Omega$ is symmetric or skew-symmetric;
\item[---] the automorphism $\sigma$ is an involution
and there exists $a \in F\setminus \{0\}$ such that
the form $a\Omega$ is hermitian.
\end{enumerate}
$\blacksquare$
}\end{rem}

\subsection{Symplectic  forms}

A bilinear form $\Omega$ on $V$
will be called {\it symplectic} if
$$
\Omega(x,x)=0\;\;\;\forall\; x\in V\;.
$$
This condition implies that
$$\Omega(x+y,x+y)=\Omega(y,x)+\Omega(x,y)=0 \;\;\;\forall\; x,y
\in V\;.$$
Therefore, any symplectic form is skew-symmetric.

The inverse statement does not hold.
For a skew-symmetric form
$\Omega$ on $V$ we have the
equality
$$\Omega(x,x)=-\Omega(x,x)\;\;\;\forall\;x\in V$$
showing that $\Omega$
is symplectic if the characteristic of
the field $F$
is not equal to $2$.
In the case when $char F=2$ there exists a skew-symmetric
form which is not symplectic
(see Example 1.1.10).

\begin{exmp}{\rm
On each even-dimensional vector space  there exists
a symplectic form.
For a $2k$-dimensional vector space
fix some base $\{x_{i},y_{i}\}^{k}_{i=1}$.
Then the $(Id,Id)$-bilinear form $\Omega$ on $V$ defined
by the conditions
$$\Omega(x_{i},x_{j})=\Omega(y_{i},y_{j})=0
\;\;\;\forall\; i,j=1,...,k\;,$$
$$\Omega(x_{i},y_{j})=0\mbox{ if }i\ne j\;,$$
$$\Omega(x_{i},y_{i})=1\;\;\;\forall\;i=1,...,k$$
is symplectic.
}$\blacksquare$
\end{exmp}

In what follows for a symplectic form $\Omega$
defined on some $2k$-dimensional vector space
each base $\{x_{i},y_{i}\}^{k}_{i=1}$
satisfying the conditions from Example 1.1.9
will be called an $\Omega$-{\it base}.

\begin{prop}
On odd-dimensional vector spaces
there are not $($non-singular$)$
symplectic bilinear forms.
For each
symplectic form $\Omega$ defined on
an even-dimensional vector space
there exists an $\Omega$-base.
\end{prop}

{\bf Proof.}
For one-dimensional
and two-dimensional vector spaces
our statement is trivial.
Assume that it holds for
each vector space the dimensional of
which is less than $n$
and consider a non-singular symplectic form
$\Omega$ defined on an $n$-dimensional space
$V$.
For a vector $x\in V\setminus\{0\}$
there exists a vector $y\in V$
satisfying the condition
$\Omega(x,y)=1$. Denote by $U$
the subspace generated by $x$ and $y$.
Then
$$U\cap U^{\perp}_{\Omega}=\{0\}$$
and the restriction of $\Omega$
on $U^{\perp}_{\Omega}$ is
a non-singular symplectic form on
$U^{\perp}_{\Omega}$;
denote it by $\Omega'$.
The equality $\dim U^{\perp}_{\Omega}=n-2$
and the inductive hypothesis
guarantee that the number $n$
is even and there exists an $\Omega'$-base
$B$.
Then $\{x,y\}\cup B$ is an
$\Omega$-base.
$\blacksquare$

\begin{cor}
Any two symplectic forms $\Omega'$
and $\Omega''$
defined
on some even-dimensional vector space $V$
are similar.
\end{cor}

{\bf Proof.}
For each $f\in \LL(V)$
transferring an $\Omega''$-base
to an $\Omega'$-base we have
the equality $f^{*}(\Omega')=\Omega''$.
$\blacksquare$

\begin{exmp}{\rm
Assume that
$$char F=2 \mbox{ and }\dim V =2k+1\;.$$
Let $\{x_{i}\}^{2k+1}_{i=1}$ be a base for $V$.
Denote by $V'$  the $2k$-dimensional linear subspace
generated by the vectors $x_{2}$,...,$x_{2k+1}$
and consider a symplectic form $\Omega'$ defined  on
$V'$ (see Example 1.1.9).
The bilinear form $\Omega$ on $V$
defined by the conditions
$$\Omega(x_{1},x)=0\;\;\;\forall\;x\in
V'\;,$$
$$\Omega(x,y)=\Omega'(x,y)\;\;\;\forall\;x,y\in V'\;,$$
$$\Omega(x_{1},x_{1})=1=-1$$
is skew-symmetric.
However, it is not symplectic.
$\blacksquare$
}\end{exmp}

\section{Three classes of transformations of the Grassmannian
manifolds}
\setcounter{equation}{0}
\setcounter{prop}{0}
\setcounter{lemma}{0}
\setcounter{rem}{0}
\setcounter{exmp}{0}
\setcounter{defn}{0}

\subsection{Transformations of the Grassmannian manifolds
defined by linear transformations}

Recall that the  {\it Grassmannian manifold}
$\ggg^{n}_{k}(V)$ is a set of $k$-dimensional linear
subspaces of an $n$-dimensional vector space $V$.
In what follows a linear subspace of $V$
(an element of the respective Grassmannian manifold)
will be called a plane.
If its dimension is equal to $1$ then
we say that it is a line.

Any  linear isomorphism $f \in \LL(V,V')$
defines the bijection
$L^{n}_{k}(f)$ of $\ggg^{n}_{k}(V)$ onto $\ggg^{n}_{k}(V')$
which will be called a {\it linear isomorphism}
of $\ggg^{n}_{k}(V)$ onto $\ggg^{n}_{k}(V')$.
Denote by $\LL^{n}_{k}(V,V')$ the class of
linear isomorphisms of $\ggg^{n}_{k}(V)$ onto
$\ggg^{n}_{k}(V')$.

If $V=V'$ then we write $\LL^{n}_{k}(V)$.
Elemensts of $\LL^{n}_{k}(V)$
are called {\it linear transformations}
of $\ggg^{n}_{k}(V)$.

The equality
$L^{n}_{k}(f_{1})=L^{n}_{k}(f_{2})$
holds for some linear isomorphisms $f_{1},f_{2}\in \LL(V,V')$
if and only if there exists $a \in F\setminus \{0\}$
such that $f_{2}=af_{1}$.

For
$f \in \LL^{n}_{k}(V,V')$ consider
a linear isomorphism $f_{1}\in \LL(V,V')$
such that $f=L^{n}_{k}(f_{1})$
and define
$$L^{n}_{k\,m}(f)=L^{n}_{m}(f_{1})\;.$$
For other linear
isomorphism $f_{2}\in \LL(V,V')$
satisfying the condition
$f=L^{n}_{k}(f_{2})$
there exists $a\in F\setminus\{0\}$
such that $f_{2}=af_{2}$. Then
$L^{k}_{m}(f_{1})=L^{k}_{m}(f_{2})$.
Therefore,
the mapping
$L^{n}_{k\,m}$
is well-defined.

\begin{prop}
The mapping $L^{n}_{k\,m}$ is a bijection
of $\LL^{n}_{k}(V,V')$ onto \\ $\LL^{n}_{m}(V,V')$.
\end{prop}

{\bf Proof.}
Let $f\in \LL^{n}_{m}(V,V')$ and
$f'\in \LL(V,V')$ be a linear isomorphism
satisfying the condition $L^{n}_{m}(f')=f$.
Then
$$L^{n}_{k\,m}(L^{n}_{k}(f'))=f\;;$$
i.e. we have proved that $L^{n}_{k\,m}$
is a mapping onto.

Now assume that the equality
$$L^{n}_{k\,m}(f_{1})=L^{n}_{k\,m}(f_{2})$$
holds for some $f_{1},f_{2}\in \LL^{n}_{k}(V,V')$
and consider $f'_{1},f'_{2}\in \LL(V,V')$
such that
$$L^{n}_{k}(f'_{i})=f_{i}\;\;\;i=1,2\;.$$
Then $L^{n}_{m}(f'_{1})=L^{n}_{m}(f'_{1})$
and there exists $a\in F\setminus\{0\}$
such that $f'_{2}=af'_{1}$.
Therefore, $f_{1}=f_{2}$.
$\blacksquare$

The similar arguments show that
$$(L^{n}_{k\,m})^{-1}=L^{n}_{m\,k}$$
and
$$L^{n}_{m\,m'}L^{n}_{k\,m}=L^{n}_{k\,m'}\;.$$
We have also
$$ (L^{n}_{k}(f))^{-1}=L^{n}_{k}(f^{-1})\;\;\;
\forall\;f\in \LL(V)\;,
$$
$$
L^{n}_{k}(f_{2})L^{n}_{k}(f_{1})=
L^{n}_{k}(f_{2}f_{1})\;\;\;
\forall\;f_{1},f_{2}\in \LL(V)\;,
$$
$$ (L^{n}_{k\,m}(f))^{-1}=L^{n}_{k\,m}(f^{-1})\;\;\;
\forall\;f\in \LL^{n}_{k}(V)\;,
$$
$$
L^{n}_{k\,m}(f_{2})L^{n}_{k\,m}(f_{1})=
L^{n}_{k\,m}(f_{2}f_{1})\;\;\;
\forall\;f_{1},f_{2}\in \LL^{n}_{k}(V)\;.
$$
In other words, the following statement
is fulfilled.

\begin{prop}
The composition operation
defines on $\LL^{n}_{k}(V)$ the group
structure such that
$L^{n}_{k}$ is  a homomorphism
of the group $\LL(V)$
onto the group $\LL^{n}_{k}(V)$,
$$Ker L^{n}_{k}=\{\;aId\;|\;a\in F\setminus\{0\}\;\}$$
and $L^{n}_{k\,m}$ is an isomorphism
of $\LL^{n}_{k}(V)$ onto $\LL^{n}_{m}(V)$.
\end{prop}

{\bf Proof.}
It is a trivial consequence of the last
four equalities.
$\blacksquare$

\begin{rem}{\rm
In Remark 1.1.1 we considered  the isomorphism
of $\LL(V)$ onto $\LL(V')$ defined by  some
linear isomorphism $f \in \LL(V,V')$.
It transfers any $aId_{V}$ ($a \in F$)
to $aId_{V'}$. This implies that
$L^{n}_{k}(f)$ defines the similar
isomorphism of $\LL^{n}_{k}(V)$ onto
$\LL^{n}_{k}(V')$.
$\blacksquare$
}\end{rem}

Let $s\in \ggg^{n}_{m}(V)$. Define
$$\ggg^{n}_{k}(s)=\left \{ \begin{array}{ll} \{\;l
\in \ggg^{n}_{k}(V)\;|\; l \subset s \;\}\;\;\; \mbox{ if }\;
m >k
\\ \{\;l\in \ggg^{n}_{k}(V)\;|\; s \subset l \;\}\;\;\;
\mbox{ if }\;m <k\;.
\end{array}
\right.$$
It is easy to see that for any two
planes $l\in \ggg^{n}_{k}(V)$ and
$s\in \ggg^{n}_{m}(V)$ the conditions
$l\in \ggg^{n}_{k}(s)$ and $s\in \ggg^{n}_{m}(l)$
are equivalent.

\begin{lemma}
For each linear isomorphism
$f\in \LL^{n}_{k}(V,V')$ the following equality
$$
f(\ggg^{n}_{k}(s))=\ggg^{n}_{k}(L^{n}_{k\,m}(f)(s))
\;\;\;\forall\;s\in \ggg^{n}_{m}(V)
$$
holds true.
\end{lemma}

{\bf Proof.}
It is trivial.
$\blacksquare$

\subsection{Bijections of
$\ggg^{n}_{k}(V)$ onto $\ggg^{n}_{n-k}(V)$
defined by bilinear forms}

A bilinear form $\Omega$ on $V$ defines the bijection
$F^{n}_{k\,n-k}(\Omega)$
of $\ggg^{n}_{k}(V)$ onto $\ggg^{n}_{n-k}(V)$
tranferring each plane $s\in \ggg^{n}_{k}(V)$ to the
$\Omega$-orthogonal complement
$s^{\perp}_{\Omega}\in \ggg^{n}_{n-k}(V)$.
In what follows we say that it is a {\it bijection of
$\ggg^{n}_{k}(V)$ onto $\ggg^{n}_{n-k}(V)$ defined by
the bilinear form $\Omega$.}

Let us consider some trivial properties of it.
\begin{enumerate}
\item[1.] The
equality
$$F^{n}_{k\,n-k}(\Omega)=F^{n}_{k\,n-k}(a\Omega)$$
holds for each $a\in F\setminus\{0\}$.
\item[2.] Let
$\Omega'$ be the inverse bilinear form. Then
$$(F^{n}_{k\,n-k}(\Omega))^{-1}=
F^{n}_{n-k\,k}(\Omega')\;.$$
\item[3.] The equality
$$(F^{n}_{k\,n-k}(\Omega))^{-1}=F^{n}_{n-k\,k}(\Omega)$$
are fulfilled
if and only if the form $\Omega$ is
reflexive.
\end{enumerate}

For the mapping $F^{n}_{k\,n-k}(\Omega)$
we have the next analogy of Lemma 1.2.1.

\begin{lemma}
The equality
$$F^{n}_{k\,n-k}(\Omega)(\ggg^{n}_{k}(s))=
\ggg^{n}_{n-k}(F^{n}_{m\,n-m}(\Omega)(s))$$
holds for any plane $s\in \ggg^{n}_{m}(V)$.
\end{lemma}

{\bf Proof.}
It is trivial.
$\blacksquare$

Lemmas 1.2.1 and 1.2.2 show that
in the case $n=2k$
the transformation
$F^{2k}_{k\,k}(\Omega)$ of
$\ggg^{2k}_{k}(V)$
is not linear.

The mapping
$F^{n}_{k\,n-k}(\Omega)$ could be described in other terms.
Denote by $\Phi^{n}_{k\,n-k}$
the bijection of $\ggg^{n}_{k}(V)$ onto $\ggg^{n}_{n-k}(V^{*})$
transferring each plane $s\in \ggg^{n}_{k}(V)$ to the orthogonal
complement $s^{\perp}\in \ggg^{n}_{n-k}(V^{*})$.
The similar bijection of
$\ggg^{n}_{k}(V^{*})$ onto $\ggg^{n}_{n-k}(V^{**}=V)$
will be denoted by
$(\Phi^{*})^{n}_{k\,n-k}$.
Then
$$
(\Phi^{*})^{n}_{n-k\,k}\Phi^{n}_{k\,n-k}=Id
$$
and
$$
\Phi^{n}_{n-k\,k}(\Phi^{*})^{n}_{k\,n-k}=Id\;.
$$
An immediate verification shows that $F^{n}_{k\,n-k}(\Omega)$
could be represented as the composition
$$F^{n}_{k\,n-k}(\Omega)=(\Phi^{*})^{n}_{k\,n-k}L^{n}_{k}(\Omega_{1})=
L^{n}_{n-k}(\Omega^{-1}_{2})\Phi^{n}_{k\,n-k}$$
(the linear isomorphisms
$\Omega_{1}$ and $\Omega_{2}$
were defined in Subsection 1.1.3).
Inversely, for any two linear isomorphisms
$f \in \LL(V,V^{*})$ and $g \in \LL(V^{*},V)$
there exist bilinear forms $\Omega'$
and $\Omega''$ such that
$\Omega'_{1}=f$ and $\Omega''_{2}=g^{-1}$.
Then
$$(\Phi^{*})^{n}_{k\,n-k}L^{n}_{k}(f)=F^{n}_{k\,n-k}(\Omega')\;,$$
$$L^{n}_{n-k}(g)\Phi^{n}_{k\,n-k}=F^{n}_{k\,n-k}(\Omega'')\;.$$

Let $n=2k$. Denote by $\FF^{2k}_{k}(V)$
the union of $\LL^{2k}_{k}(V)$ and the class of
transformations of $\ggg^{2k}_{k}(V)$
defined by bilinear forms on $V$.
Then the following statement holds true.

\begin{prop}
The composition operation  defines on
$\FF^{2k}_{k}(V)$ the group structure
such that
$\LL^{2k}_{k}(V)$ is a normal subgroup
of $\FF^{2k}_{k}(V)$ and
$$\FF^{2k}_{k}(V) / \LL^{2k}_{k}(V)\cong \zz_{2}\;.$$
\end{prop}

{\bf Proof.}
Let $\Omega$ be a bilinear form on $V$
and $f\in \LL^{2k}_{k}(V)$.
Consider $F^{2k}_{k\,k}(\Omega)$
as the composition of a linear
isomorphism belonging to $\LL^{2k}_{k}(V,V^{*})$
and $(\Phi^{*})^{2k}_{k\,k}$.
Then $F^{2k}_{k\,k}(\Omega)f$ is  the composition of a linear
isomorphism belonging to $\LL^{2k}_{k}(V,V^{*})$
and $(\Phi^{*})^{2k}_{k\,k}$;
i.e. it is defined by some bilinear form on $V$.
The mapping $F^{2k}_{k\,k}(\Omega)$
could be represented as the composition of
$\Phi^{2k}_{k\,k}$ and a linear isomorphism
belonging to $\LL^{2k}_{k}(V^{*},V)$.
This implies that $fF^{2k}_{k\,k}(\Omega)$
has the similar representation;
therefore, it is
defined by
a bilinear form.

For any two bilinear forms
$\Omega'$ and $\Omega''$ on $V$ there exist
$f\in \LL^{2k}_{k}(V,V^{*})$ and $g \in
\LL^{2k}_{k}(V^{*},V)$
such that
$$F^{2k}_{k\,k}(\Omega')=(\Phi^{*})^{2k}_{k\,k}f$$
and
$$F^{2k}_{k\,k}(\Omega'')=g\Phi^{2k}_{k\,k}\;.$$
Then
$$F^{2k}_{k\,k}(\Omega'')F^{2k}_{k\,k}(\Omega')=
gf\in \LL^{2k}_{k}(V)\;.$$
We have proved that the composition of two
elements of $\FF^{2k}_{k}(V)$
belongs to $\FF^{2k}_{k}(V)$.
Properties 2 guarantees the existence of the inverse element
for each $F^{2k}_{k\,k}(\Omega)$
and first statement is proved.

Let us prove the inclusion
$$(F^{2k}_{k\,k}(\Omega))^{-1}\LL^{2k}_{k}(V)F^{2k}_{k\,k}(\Omega)
\subset \LL^{2k}_{k}(V)$$
showing that the subgroup
$\LL^{2k}_{k}(V)$ is normal.
Consider
$F^{2k}_{k\,k}(\Omega)$
as the composition
$$F^{2k}_{k\,k}(\Omega)=(\Phi^{*})^{2k}_{k\,k}f=g\Phi^{2k}_{k\,k}\;,$$
where $f\in \LL^{2k}_{k}(V,V^{*})$ and
$g \in\LL^{2k}_{k}(V^{*},V)$.
Then
$$(F^{2k}_{k\,k}(\Omega))^{-1}\LL^{2k}_{k}(V)F^{2k}_{k\,k}(\Omega)=
f^{-1}\Phi^{2k}_{k\,k}\LL^{2k}_{k}(V)g\Phi^{2k}_{k\,k}\subset$$
$$f^{-1}\Phi^{2k}_{k\,k}\LL^{2k}_{k}(V^{*},V)\Phi^{2k}_{k\,k}
\subset
f^{-1}\Phi^{2k}_{k\,k}(\Phi^{*})^{2k}_{k\,k}\LL^{2k}_{k}(V,V^{*})
\subset$$
$$f^{-1}\LL^{2k}_{k}(V,V^{*})\subset \LL^{2k}_{k}(V)\;.$$

It was proved above that
for any two  bilinear forms
$\Omega'$ and $\Omega''$ on $V$
the composition
$$F^{2k}_{k\,k}(\Omega'')F^{2k}_{k\,k}(\Omega')$$
is a linear transformation of $\ggg^{2k}_{k}(V)$.
Therefore,
the factor group
$$\FF^{2k}_{k}(V) / \LL^{2k}_{k}(V)$$
contains only two classes.  One of them coincides with
$\LL^{2k}_{k}(V)$ and second class consists of transformations
of $\ggg^{2k}_{k}(V)$
defined by bilinear forms.
$\blacksquare$

\subsection{Transformations of $\ggg^{n}_{m}(V)$
induced by  transformations of $\ggg^{n}_{k}(V)$}

Let $f$ be a transformation
of $\ggg^{n}_{k}(V)$ and for any plane
$s \in\ggg^{n}_{m}(V)$ ($m\ne k$)
there exists a plane
$s_{f} \in\ggg^{n}_{m}(V)$
such that
$$f(\ggg^{n}_{k}(s))=\ggg^{n}_{k}(s_{f})\;.$$
Denote by $g$ the mapping transferring each
$s\in\ggg^{n}_{m}(V)$ to $s_{f}$.
It is trivial that $g$ is a bijection
of $\ggg^{n}_{m}(V)$ into itself
(otherwise, the mapping $f$ is not bijective).
In that follows we restrict ourself only to the case
when $g$ is a transformation  of
$\ggg^{n}_{m}(V)$;
i.e. for any plane $s \in\ggg^{n}_{m}(V)$
there exists a plane
$s'_{f} \in\ggg^{n}_{m}(V)$
such that
$$f^{-1}(\ggg^{n}_{k}(s))=\ggg^{n}_{k}(s'_{f})\;.$$
In this case we say that the transformation $f$ {\it induces}
the transformation $g$
or $g$ is {\it induced} by $f$.

For any two planes $l \in\ggg^{n}_{k}(V)$
and $s \in\ggg^{n}_{m}(V)$ we have the following
chain of implications
$$s \in\ggg^{n}_{m}(l)\Leftrightarrow l \in\ggg^{n}_{k}(s)
\Leftrightarrow f(l) \in\ggg^{n}_{k}(g(s))
\Leftrightarrow g(s) \in\ggg^{n}_{m}(f(l))$$
showing that
$$g(\ggg^{n}_{m}(l))=\ggg^{n}_{m}(f(l))\;\;\;
\forall\;l \in\ggg^{n}_{k}(V)\;.$$
In other words, the next lemma is proved.

\begin{lemma}
If a transformation $f$
of $\ggg^{n}_{k}(V)$
induces
a transformation $g$ of $\ggg^{n}_{m}(V)$
then $g$ induces  $f$.
\end{lemma}

It must be pointed out that
{\it a transformation of $\ggg^{n}_{m}(V)$
induced by some
transformation of $\ggg^{n}_{k}(V)$
is uniquely defined.}

\begin{exmp}{\rm
Lemma 1.2.1 shows that each linear transformation $f\in
\LL^{n}_{k}(V)$ induces $L^{n}_{k\,m}(f)$.} $\blacksquare$
\end{exmp}

\begin{lemma}
If a transformation $f$ of $\ggg^{n}_{k}(V)$
induces
a transformation
$g$ of $\ggg^{n}_{m}(V)$
and $g$ induces
a transformation
$h$ of $\ggg^{n}_{p}(V)$
then $h$ is induced by $f$.
\end{lemma}

{\bf Proof.}
We have to consider the following six cases:
\begin{enumerate}
\item[(i)] $k>m>p$,
\item[(ii)] $k>p>m$,
\item[(iii)] $m>k>p$,
\item[(iv)] $p>m>k$,
\item[(v)] $p>k>m$,
\item[(vi)] $m>p>k$.
\end{enumerate}
We restrict ourself to the cases when $k>p$
((i) -- (iii)).
Lemma 1.2.3 shows that cases (iv) -- (vi)
($p>k$)
could be reduced to them
by  tranponsing of $f$ and $h$.

Let $s \in\ggg^{n}_{p}(V)$.
In case (i) a plane $l\in \ggg^{n}_{k}(V)$
belongs to the set $\ggg^{n}_{k}(s)$ if and only if
it contains some plane $t\in \ggg^{n}_{m}(s)$.
In other words,
$$\ggg^{n}_{k}(s)=
\bigcup_{t\in \ggg^{n}_{m}(s)}\ggg^{n}_{k}(t)\;.$$
For each plane $t\in \ggg^{n}_{m}(V)$ the mapping
$f$ transfers $\ggg^{n}_{k}(t)$ to $\ggg^{n}_{k}(g(t))$.
This implies that
$$f(\ggg^{n}_{k}(s))
=\bigcup_{t\in g(\ggg^{n}_{m}(s))}\ggg^{n}_{k}(t)
=\bigcup_{t\in \ggg^{n}_{m}(h(s))}\ggg^{n}_{k}(t)
=\ggg^{n}_{k}(h(s))\;.$$
In case (ii) a plane $l\in \ggg^{n}_{k}(V)$
belongs to the set $\ggg^{n}_{k}(s)$ if and only if
it contains each plane belonging to $\ggg^{n}_{m}(s)$.
We obtain
the equality
$$\ggg^{n}_{k}(s)=\bigcap_{t\in \ggg^{n}_{m}(s)}
\ggg^{n}_{k}(t)$$
showing that
$$f(\ggg^{n}_{k}(s))=\ggg^{n}_{k}(h(s))\;.$$
For cases (i) and (ii)
our statement is proved.

Consider case (iii).
Let $\Omega$ be a bilinear form on $V$ and
$$f'=F^{n}_{k\,n-k}(\Omega)f(F^{n}_{k\,n-k}(\Omega))^{-1}\;,$$
$$g'=F^{n}_{m\,n-m}(\Omega)g(F^{n}_{m\,n-m}(\Omega))^{-1}\;,$$
$$h'=F^{n}_{p\,n-p}(\Omega)h(F^{n}_{p\,n-p}(\Omega))^{-1}\;.$$
An immediate verification shows that
$f'$ induces $g'$ and $g'$ induces $h'$.
Moreover, $f'$ induces $h'$
if and only if $f$ induces $h$.
This implies that case (iii)
could be reduced to case (v).
Recall that the last case was reduced to
case (ii).
$\blacksquare$

The next theorem will be proved in Section 1.4.
To prove it we exploit the
Fundamental Theorem of Projective Geometry
and the Chow Theorem.

\begin{theorem}
If a transformation $f$ of $\ggg^{n}_{k}(V)$
induces a transformation
of $\ggg^{n}_{m}(V)$ and $m\ne k$
then $f\in \LL^{n}_{k}(V)$.
\end{theorem}

\section{Fundamental Theorem of Projective \\Geometry}
\setcounter{equation}{0}
\setcounter{prop}{0}
\setcounter{lemma}{0}
\setcounter{rem}{0}
\setcounter{theorem}{0}
\setcounter{exmp}{0}
\setcounter{defn}{0}

\subsection{Formulation and remarks}

Let $x \in V\setminus\{0\}$.  Denote by $l(x)$ the line containing
the vector $x$.
Lines
$$l_{1},...,l_{k}\in \ggg^{n}_{1}(V)$$
will be called {\it linearly independent} if
there exists a collection of linearly independent  vectors
$x_{1},...,x_{k}$ such that
$$l_{i}=l(x_{i})\;\;\;\forall\;i=1,...,k\;.$$
Otherwise, we say that the lines
are {\it linearly dependent}.

A linear transformation of $\ggg^{n}_{1}(V)$
maps each collection of
linearly independent
lines to a collection of
linearly independent
lines.
The inverse statement is known as the
Fundamental Theorem of Projective Geometry.

\begin{theorem}
{\rm (the Fundamental Theorem of Projective Geometry)}
Let $\dim V \ge 3$ and $f$ be a transformation
of $\ggg^{n}_{1}(V)$
such that
$f$ and $f^{-1}$
map each collection of
linearly independent
lines to a collection
of linearly independent
lines.
Then $f\in \LL^{n}_{1}(V)$.
\end{theorem}

\begin{rem}{\rm
For the case when $\dim V =2$ the analogous statement fails.
In this case any three lines are linearly dependent
and each transformation
of $\ggg^{2}_{1}(V)$
maps any collection of
linearly independent
lines to a collection
of linearly independent
lines.
$\blacksquare$
}\end{rem}

The next lemma will be used in what follows to reformulate
Theorem 1.3.1 in other terms.

\begin{lemma}
For a transformation $f$ of $\ggg^{n}_{1}(V)$
the following
conditions are equivalent:
\begin{enumerate}
\item[$(a)$] for any plane
$s\in \ggg^{n}_{k}(V)$
there exist planes
$s_{f},s'_{f}\in \ggg^{n}_{k}(V)$
such that
$$
f(\ggg^{n}_{1}(s))=\ggg^{n}_{1}(s_{f})$$
and
$$f^{-1}(\ggg^{n}_{1}(s))=\ggg^{n}_{1}(s'_{f})\;;$$
\item[$(b)$] the mappings $f$ and $f^{-1}$
transfer each collection
of $m$ linearly independent
lines $(m\le k+1)$ to a collection
of linearly independent
lines.
\end{enumerate}
\end{lemma}

{\bf Proof.}
$(a)\Rightarrow (b)$.
Condition (a) implies that
$f$ and $f^{-1}$ map any collection of
$k+1$ linearly dependent lines  to a collection of
linearly dependent lines,
since $k+1$ lines are linearly dependent if and only if
they are contained in some $k$-dimensional plane.
Therefore, our mappings transfer any collection of
$k+1$ linearly independent lines  to a collection of
linearly independent lines.

For linearly independent lines
$$l_{1},...,l_{m}\in \ggg^{n}_{1}(V)\;(m\le k)$$
there exists lines
$$l_{m+1},...,l_{k+1}\in \ggg^{n}_{1}(V)$$
such that $l_{1},...,l_{k+1}$ is a collection
of linearly independent lines.
Then the lines $f(l_{1}),...,f(l_{k+1})$
are linearly independent. This implies that
$f(l_{1}),...,$ $f(l_{m})$ are
linearly independentent.
It is trivial that the analogous statement
holds for the mapping $f^{-1}$.

$(b)\Rightarrow (a)$.
Let $l_{1},...,l_{k}$ be lines
generating a plane $s \in \ggg^{n}_{k}(V)$.
Condition $(b)$ guarantees that the
lines $f(l_{1}),...,f(l_{k})$
generate some plane $s_{f} \in \ggg^{n}_{k}(V)$.
If $l\in \ggg^{n}_{1}(s)$ then
$$l_{1},...,l_{k},l \mbox{ and }f(l_{1}),...,f(l_{k}),f(l)$$
are
collections of linearly dependent lines. Therefore, $f(l)\in
\ggg^{n}_{1}(s_{f})$ and
we obtain the inclusion
$$f(\ggg^{n}_{1}(s))\subset\ggg^{n}_{1}(s_{f})\;.$$
The proof of the inverse
inclusion is similar.
These arguments imply also the existence of a plane $s'_{f}\in
\ggg^{n}_{k}(V)$ satisfying second condition.
$\blacksquare$

Now Theorem 1.3.1 could be reformulated in the following form.
It could be considered as a special case of Theorem
1.2.1.

\begin{theorem}
Let $\dim V \ge 3$ and $f$ be a transformation of
$\ggg^{n}_{1}(V)$ inducing some  transformation of
$\ggg^{n}_{n-1}(V)$.
Then $f\in \LL^{n}_{1}(V)$.
\end{theorem}

{\bf Proof.}
A transformation $f$ of $\ggg^{n}_{1}(V)$
satisfies  condition (a)
from Lemma 1.3.1 if and only if it induces
a  transformation of $\ggg^{n}_{k}(V)$.
Therefore, the required statement is
a consequence of Theorem 1.3.1
and Lemma 1.3.1.
$\blacksquare$

\begin{rem}{\rm
Theorem 1.3.2 and Lemma 1.2.3 show that
if a transformation
$f$ of $\ggg^{n}_{n-1}(V)$
($\dim V \ge 3$) induces a transformation of $\ggg^{n}_{1}(V)$
then $f\in \LL^{n}_{n-1}(V)$.
$\blacksquare$
}\end{rem}

\subsection{Proof of Theorem 1.3.1}

Our proof will be decomposed on a few steps.

{\it Fist step.}
Fix a base $\{x_{i}\}^{n}_{i=1}$ for $V$
and for each $i=1,...,n$ define $l_{i}=l(x_{i})$.
Then the lines
$$l'_{1}=f(l_{1}),...,l'_{n}=f(l_{n})$$
are linearly independent and there exists a base
$\{y'_{i}\}^{n}_{i=1}$ for $V$ such that
$l'_{i}=l(y'_{i})$ for any $i=1,...,n$.

The mapping $f$ transfers a line contained in the plane
generated by $x_{1}$ and $x_{i}$ to
a line contained in the plane
generated by $y'_{1}$ and $y'_{i}$
(Lemma 1.3.1).
This implies the existence of
$a_{2},...,a_{n} \in F$ such
that
$$f(l(x_{1}+x_{i}))=l(y'_{1}+a_{i}y'_{i})\;\;\;
\forall\;i=2,...,n\;.$$
Define
$$y_{1}=y'_{1} \mbox{ and }
y_{i}=a_{i}y'_{i}\;\;\;\forall\;i=2,...,n\;.$$
Then $\{y_{i}\}^{n}_{i=1}$ is a base and
$l'_{i}=l(y_{i})$ for each $i=1,...,n$.
Moreover,
$$f(l(x_{1}+x_{i}))=l(y_{1}+y_{i})\;\;\;
\forall\;i=2,...,n\;.$$

In what follows we denote by
$l+s$ the two-dimensional plane
generated by lines $l$ and $s$.
In these terms we have the following
inclusion
$$ l(x_{i}-x_{j}) \subset \left \{ \begin{array}{l}
l_{i}+l_{j}
\\ l(x_{1}+x_{i})+l(x_{1}+x_{j})
\end{array}
\right.$$
showing that
$$
f(l(x_{i}-x_{j})) \subset \left \{ \begin{array}{l}
l'_{i}+l'_{j}
\\ l(y_{1}+y_{i})+l(y_{1}+y_{j})\;.
\end{array}
\right.$$
The last inclusion guarantees that
for a vector $y\in V$ satisfying the condition
$$l(y)=f(l(x_{i}-x_{j}))$$
there exist $a,b\in F$
such that
$$y=a(y_{1}+y_{i})+b(y_{1}+y_{j})\;.$$
However, $y$ is a linear combination
of the vectors $y_{i}$ and $y_{j}$.
Therefore, $a=-b$ and
\begin{equation}
f(l(x_{i}-x_{j}))=l(y_{i}-y_{j})\;\;\;\forall\;i,j=2,...,n\;.
\end{equation}

{\it Second step.}
Lemma 1.3.1 shows that
for each
$a \in F$ and $i=2,...,n$ there exists
$\sigma_{i}(a)\in F$
such that
$$f(l(x_{1}+ ay_{i}))=l(y_{1}+\sigma_{i}(a)y_{i})\;.$$
This equality
defines the transformation $\sigma_{i}$ of $F$
for each $i=2,...,n$.
Prove that $\sigma_{i}=\sigma_{j}$
for any two numbers $i$ and $j$.

Consider $a \in F$ and define
$a_{i}=\sigma_{i}(a)$, $a_{j}=\sigma_{j}(a)$.
The  inclusion
$$
l(ax_{i}-ax_{j}) \subset \left \{ \begin{array}{l}
l_{i}+l_{j}
\\ l(x_{1}+ax_{i})+l(x_{1}+ax_{j})\;.
\end{array}
\right.$$
implies that
$$
f(l(ax_{i}-ax_{j})) \subset \left \{ \begin{array}{l}
l'_{i}+l'_{j}
\\ l(y_{1}+a_{i}x_{i})+l(y_{1}+a_{j}y_{j})\;.
\end{array}
\right.$$
Then
$$f(l(x_{i}-x_{j}))=f(l(ax_{i}-ax_{j}))=
l(a_{i}y_{i}-a_{j}y_{j})$$
(the proof of this equality is similar to
the proof of equation (1.3.1)).
Equation (1.3.1) shows that
$$l(y_{i}-y_{j})=l(a_{i}y_{i}-a_{j}y_{j})\;.$$
Therefore,
$a_{i}=a_{j}$.
In what follows we write $\sigma$ in place of each
$\sigma_{i}$.

{\it Third step.}
Let
$x=x_{1}+\Sigma^{n}_{j=2} a_{j}x_{j}$.
Then
$$
l(x)\subset l(x_{1}+a_{i}x_{i})+l(x')\;,
$$
$$
f(l(x))\subset l(y_{1}+\sigma(a_{i})y_{i})+l(y')\;,
$$
where $x'$ and $y'$ are a linear combinations of the vectors
$$x_{2},...,x_{i-1},x_{i+1},...,x_{n}$$
and
$$y_{2},...,y_{i-1},y_{i+1},...,y_{n}\;,$$
respectively.
This implies that
for a vector $y=y_{1}+\Sigma^{n}_{j=2} b_{j}y_{j}$
satisfying the condition $l(y)=f(l(x))$
we have
$b_{i}=\sigma(a_{i})$
for each $i=2,...,n$.
In other words,
$$f(l(x_{1}+\Sigma^{n}_{j=2} a_{j}x_{j}))=
l(y_{1}+\Sigma^{n}_{j=2} \sigma(a_{j})y_{j})\;.$$

Let us consider the case when
$x=\Sigma^{n}_{j=2} a_{j}x_{j}$.
It is trivial that
$$l(x)\subset l(x_{1}+\Sigma^{n}_{j=2} a_{j}x_{j})+l_{1}\;,$$
$$f(l(x))\subset l(y_{1}+\Sigma^{n}_{j=2}
\sigma(a_{j})y_{j})+l'_{1}\;.$$
The last inclusion shows that
for a vector $y\in V$ satisfying the condition
$l(y)=f(l(x))$ there exist $a,b\in F$
such that
$$y=a(y_{1}+\Sigma^{n}_{j=2}\sigma(a_{j})y_{j})
+by_{1}\;.$$
The line $f(l(x))$ is contained in the plane
generated by the vectors $y_{2},...,y_{n}$.
Therefore, $a=-b$ and
$$f(l(\Sigma^{n}_{j=2} a_{j}x_{j}))=
l(\Sigma^{n}_{j=2} \sigma(a_{j})y_{j})\;.$$

{\it Fourt step.}
Now we can prove that $\sigma \in Aut(F)$.
It is trivial that  $\sigma(0)=0$ and $\sigma(1)=1$.
For any two $a,b \in F$ we have
$$l(x_{1}+ (a+b)x_{2}+x_{3})\subset
l(x_{1}+ax_{2})+l(bx_{2}+x_{3})\;.$$
Then
$$l(y_{1}+ \sigma(a+b)y_{2}+y_{3})=f(l(x_{1}+
(a+b)x_{2}+x_{3}))\subset$$
$$l(y_{1}+\sigma(a)y_{2})+l(\sigma(b)y_{2}+y_{3})\;.$$
This implies the existence of $c,d\in F$
such that
$$y_{1}+ \sigma(a+b)y_{2}+y_{3}=
c(y_{1}+\sigma(a)y_{2})+d(\sigma(b)y_{2}+y_{3})\;.$$
It is trivial that $c=d=1$
and
$$\sigma(a+b)=\sigma(a)+\sigma(b)\;.$$

We have also the following inclusion
$$l(x_{1}+ abx_{2}+bx_{3})\subset
l_{1}+l(abx_{2}+bx_{3})=l_{1}+l(ax_{2}+x_{3})$$
showing that
$$l(y_{1}+\sigma(ab)y_{2}+\sigma(b)y_{3})=
f(l(x_{1}+abx_{2}+bx_{3}))\subset$$
$$l'_{1}+l(\sigma(a)y_{2}+y_{3})=
l'_{1}+l(\sigma(a)\sigma(b)y_{2}+\sigma(b)y_{3})\;.$$
It is not difficult to see that
$$\sigma(ab)=\sigma(a)\sigma(b)$$
and the required statement is proved.

{\it Last step.}
If $a_{1}\ne 0$ then
$$f(l(\Sigma^{n}_{i=1}a_{i}x_{i}))=
f(l(x_{1}+\Sigma^{n}_{i=2}a_{i}a^{-1}_{1}x_{i}))=$$
$$l(y_{1}+\Sigma^{n}_{i=2}\sigma(a_{i}a^{-1}_{1})y_{i})=
l(\Sigma^{n}_{i=1}\sigma(a_{i})y_{i})$$
(for the case $a_{1}=0$ this
equality was proved above).
Then
$$f=L^{n}_{1}(gf_{\sigma B})\;,$$
where
$B=\{x_{i}\}^{n}_{i=1}$ and
$g$ is the element of $\in \LL_{Id}(V)$
transferring the base $\{x_{i}\}^{n}_{i=1}$
to the base
$\{y_{i}\}^{n}_{i=1}$ (see Example 1.1.3).
$\blacksquare$

\begin{rem}{\rm
This proof was taken from [O'M] (see also [Art]).
In [Die] there is other proof but it is more complicated.
$\blacksquare$
}\end{rem}

\section{Transformations of the Grassmannian \\ manifolds
preserving the distance \\between planes}
\setcounter{equation}{0}
\setcounter{prop}{0}
\setcounter{lemma}{0}
\setcounter{rem}{0}
\setcounter{theorem}{0}
\setcounter{exmp}{0}
\setcounter{defn}{0}

\subsection{Chow's Theorem}

In this section we give the analogy of Theorem 1.3.1 for
transformations of $\ggg^{n}_{k}(V)$,
where $1<k<n-1$. It is known as the Chow Theorem.
To formulate it we exploit the notion of
the distance between planes.

For any two planes $l_{1}, l_{2}\in \ggg^{n}_{k}(V)$ the number
$$d(l_{1}, l_{2})=k-\dim l_{1}\cap l_{2}$$
is called the {\it distance} between them.
It is easy to see that the equality $d(l_{1}, l_{2})=i$ holds if and only
if the dimension of the plane generated by $l_{1}$ and $l_{2}$ is
equal to $k+i$.  This implies the inequality
$$d(l_{1}, l_{2}) \le
\alpha^{n}_{k}= \left \{ \begin{array}{ll}
k &\mbox{ if } k \le n-k\\
n-k &\mbox{ if }  k \ge n-k
\end{array} \right.$$
Note that the case
$d(l_{1}, l_{2})=\alpha^{n}_{k}$
could be realized.

We say that the planes $l_{1}$ and $l_{2}$ are {\it adjacent}
if $d(l_{1}, l_{2})=1$.

\begin{exmp}{\rm
For the cases $k=1,n-1$
any two planes belonging to $\ggg^{n}_{k}(V)$
are adjacent.
$\blacksquare$
}\end{exmp}

In what follows we shall use the next simple lemmas.

\begin{lemma}
The following conditions are equivalent:
\begin{enumerate}
\item[---] planes $l_{1}, l_{2} \in \ggg^{n}_{k}(V)$
are adjacent,
\item[---] the set
$\ggg^{n}_{k+1}(l_{1})\cap \ggg^{n}_{k+1}(l_{2})$
consists of unique plane {\rm(}this plane is generated
by $l_{1}$ and $l_{2}${\rm)},
\item[---] the set
$\ggg^{n}_{k-1}(l_{1})\cap \ggg^{n}_{k-1}(l_{2})$
consists of unique plane {\rm(}it is the intersection
of $l_{1}$ and $l_{2}${\rm)}.
\end{enumerate}
\end{lemma}

{\bf Proof.}
It is trivial.$\blacksquare$

\begin{lemma}
For planes $l,{\hat l} \in \ggg^{n}_{k}(V)$
the distance $d(l,{\hat l})$
is the smallest number $i$ such that
there exists a collection of
planes
$$l_{0},l_{1},...,l_{i} \in \ggg^{n}_{k}(V)\;,$$
where
$l_{0}=l$, $l_{i}={\hat l}$
and the planes $l_{j},l_{j+1}$ are
adjacent for any $j=0,...,i-1$.
\end{lemma}

{\bf Proof.}
If $d(l,{\hat l})=i$
then $l$ and ${\hat l}$ generate some
$(k+i)$-dimensional plane $s$.
There exists
a collection of linearly independent
vectors
$x_{1},...,x_{k+i}\in V$
generatying $s$ and
such that $l$ is generated by
$x_{1},...,x_{k}$ and ${\hat l}$ is generated by
$x_{i+1},...,x_{k+i}$.  For each $j=0,...i$ denote by $l_{j}$ the
plane genarated by $x_{j+1},...,x_{k+j}$. It is easy to see that
the planes $l_{0},l_{1},...,l_{i}$ satisfy the required
conditions.

Let
$$l'_{0},l'_{1},...,l'_{m} \in \ggg^{n}_{k}(V)$$
be other collections of planes
such that $l'_{0}=l$, $l'_{m}={\hat l}$
and $l'_{j},l'_{j+1}$ are
adjacent for any $j=0,...,m-1$.
Then there exist vectors
$x'_{1},...,x'_{k+m}\in V$
such that for each $j=0,...,m$
the plane $l'_{j}$ is generated
by $x'_{j+1},...,x'_{k+j}$.
This implies that
$l$ and ${\hat l}$ are contained in
the plane generated by $x'_{1},...,x'_{k+m}$.
Therefore, $d(l,{\hat l})\le m$ and
our statement is proved.
$\blacksquare$

\begin{lemma}
For a transformation
$f$ of $\ggg^{n}_{k}(V)$
the next two conditions are equivalent:
\begin{enumerate}
\item[---] $f$ and $f^{-1}$ maps
any two adjacent planes to adjacent planes;
\item[---] $f$ preserves the distance between planes, i.e.
$$d(l_{1}, l_{2})=d(f(l_{1}),f(l_{2}))\;\;\;\forall\; l_{1},
l_{2}\in \ggg^{n}_{k}(V)\;.$$ \end{enumerate} \end{lemma}

{\bf Proof.} It is a trivial consequence of Lemma 1.4.2.
$\blacksquare$

\begin{rem}{\rm
The equality $\alpha^{n}_{k}=\alpha^{n}_{n-k}$
guarantees the fulfilment of the analogous
statement for bijections  of $\ggg^{n}_{k}(V)$ onto
$\ggg^{n}_{n-k}(V)$.
$\blacksquare$}
\end{rem}

Denote by $\CC^{n}_{k}(V)$ the class of
transformations of $\ggg^{n}_{k}(V)$
preserving the distance between
planes.

\begin{exmp}{\rm
In the cases $k=1,n-1$ each transformation of
$\ggg^{n}_{k}(V)$
belongs to $\CC^{n}_{k}(V)$
(see Example 1.4.1).
$\blacksquare$
}\end{exmp}

\begin{exmp}{\rm
Let $f \in \LL^{n}_{k}(V)$.
Lemma 1.2.1 shows that
planes $l_{1},l_{2}\in \ggg^{n}_{k}(V)$
are
contained in a $(k+1)$-dimensional
plane $s$ if and only if
the planes $f(l_{1})$ and $f(l_{2})$
are contained in the plane
$L^{n}_{k\,k+1}(f)(s)$.
Then Lemma 1.4.3
guarantees the fulfilment of the inclusion
$$\LL^{n}_{k}(V)\subset \CC^{n}_{k}(V)\;.$$}
$\blacksquare$
\end{exmp}

\begin{exmp}{\rm
Now consider a bijection $g$
of $\ggg^{n}_{k}(V)$ onto $\ggg^{n}_{n-k}(V)$
defined by a bilinear form $\Omega$.
Lemma 1.2.2 states that
planes $l_{1},l_{2}\in \ggg^{n}_{k}(V)$ are
contained in some $(k+1)$-dimensional
plane $s$ if and only if
the planes
$$f(l_{1}),f(l_{2}) \in \ggg^{n}_{n-k}(V)$$
contain
the $(n-k-1)$-dimensional plane $s^{\perp}_{\Omega}$.
This implies that
$f$ preserve the distance between planes
(see Remark 1.4.1) and the following inclusion
$$\FF^{2k}_{k}(V)\subset \CC^{2k}_{k}(V)$$
holds true.
$\blacksquare$
}\end{exmp}

The Chow Theorem can be formulated
in the following form.

\begin{theorem}{\rm [Chow]}
If $1<k<n-1$ and $n\ne 2k$ then
$$\CC^{n}_{k}(V)=\LL^{n}_{k}(V)\;.$$
If $n=2k$ and $k> 1$ then
$$\CC^{2k}_{k}(V)=\FF^{2k}_{k}(V)\;.$$
\end{theorem}

We devote this section to prove the Chov theorem.
Our proof is a modification of the proof given in the Chow's
paper [Chow] (see also [Die]).

\subsection{Sets of adjacent planes}

Now we consider one special class of subsets of
$\ggg^{n}_{k}(V)$, $1<k<n-1$.  In the next subsection we use
its properties to prove the Chow Theorem.

Denote by ${\cal C}^{n}_{k}(V)$ the class of sets
$C\subset \ggg^{n}_{k}(V)$ such that any two planes
belonging to $C$ are adjacent.
Let ${\cal MC}^{n}_{k}(V)$ be the subclass  in $\CC^{n}_{k}(V)$
consists of sets $C\in {\cal C}^{n}_{k}(V)$
satisfying the next condition: any set
$C'\in {\cal C}^{n}_{k}(V)$
containing $C$ coincides with it.

\begin{prop}
For any  set $C\in {\cal C}^{n}_{k}(V)$ there exists
$C'\in {\cal MC}^{n}_{k}(V)$ containing $C$.
\end{prop}

{\bf Proof.}
Denote by ${\cal I}$ the family of sets
$C' \in {\cal C}^{n}_{k}(V)$
containing $C$.
Let ${\cal J} \subset {\cal I}$ be
a linearly ordered family (see Remark 1.4.2) and
$$C_{{\cal J}}=\bigcup_{C' \in {\cal J}}C'\;.$$
For any two planes $l_{1},l_{2}\in C_{{\cal J}}$
there exist
$C_{1},C_{2}\in {\cal J}$
containing $l_{1}$ and $l_{2}$,
respectively.
The family ${\cal J}$ is linearly ordered
and one of these set is contained in
other set.
This implies that
these planes are adjacent.
We have proved that
$$C_{{\cal J}}\in{\cal C}^{n}_{k}(V)\;;$$
i.e. the family ${\cal J}$ is
bounded above.
Then our statement is a consequence of the Zorn Lemma
(Remark 1.4.2).
$\blacksquare$

\begin{rem}{\bf Zorn's Lemma.}
{\rm
Recall the formulation of the Zorn Lem-\\ma.
Let ${\cal I}$ be a family of subsets of some set $X$
and ${\cal J}\subset {\cal I}$.
The family ${\cal J}$ is  {\it linearly ordered}
if for any two sets $U_{1}, U_{2}\in {\cal J}$
one of the inclusions
$$U_{1} \subset U_{2} \mbox{ or } U_{2} \subset U_{1}$$
holds true.
The family ${\cal J}$ is called {\it bounded above} in ${\cal I}$
if there exists $U'\in{\cal I}$
containing each $U \in {\cal J}$.
We say that a set $U\in {\cal I}$ is {\it maximal}
in the family ${\cal I}$
if any set $U'\in{\cal I}$ containing $U$
coincides with it.
The Zorn Lemma states that
{\it if any linearly ordered family in ${\cal I}$
is bounded above
then the family ${\cal I}$ has a maximal set.}
$\blacksquare$
}\end{rem}

\begin{exmp}{\rm
Let $s\in \ggg^{n}_{m}(V)$ and
$m=k-1$ or $k+1$. Then
$\ggg^{n}_{k}(s)\in {\cal C}^{n}_{k}(V)$.
Show that
$$\ggg^{n}_{k}(s)\in {\cal MC}^{n}_{k}(V)\;.$$

If $m=k+1$ then for each
plane $l$ belonging to $\ggg^{n}_{k}(V)\setminus\ggg^{n}_{k}(s)$
we have
$$\dim l\cap s \le k-1$$
and there exists a plane ${\hat l} \in \ggg^{n}_{k}(s)$
which does not contain $l\cap s$.
Then
$$\dim l\cap {\hat l} \le k-2$$
and $d(l,{\hat l})>1$. Therefore,
$$\ggg^{n}_{k}(s)\cup\{l\}\notin {\cal C}^{n}_{k}(V)
\;\;\;\forall\;l\in\ggg^{n}_{k}(V)\setminus\ggg^{n}_{k}(s)\;.$$
This implies the required.

In the case $m=k-1$ consider a bijection
$f$ of $\ggg^{n}_{k}(V)$ onto $\ggg^{n}_{n-k}(V)$
defined by a bilinear form $\Omega$. Then
$$f(\ggg^{n}_{k}(s))=\ggg^{n}_{n-k}(s^{\perp}_{\Omega})
\in {\cal MC}^{n}_{n-k}(V)\;,$$
since $s^{\perp}_{\Omega}\in \ggg^{n}_{n-k+1}(V)$.
The bijection $f$ preserves the distance between
planes and the required statement is
proved for the case $m=k-1$.
}$\blacksquare$
\end{exmp}

Now show that the class
${\cal MC}^{n}_{k}(V)$ is constituted only of sets considered in
Example 1.4.5.
In other words, the following statement holds true.

\begin{prop}
For any set $C \in {\cal MC}^{n}_{k}(V)$
there exists
a plane $s\in \ggg^{n}_{m}(V)$
such that $m=k-1$ or $k+1$ and $C=\ggg^{n}_{k}(s)$.
\end{prop}

{\bf Proof.}
Let $l_{1},l_{2}\in C$.
Define $s_{1}=l_{1} \cap l_{2}$.
Then $\dim s_{1}=k-1$, since the planes $l_{1}$ and $l_{2}$ are
adjacent.
Denote by $s_{2}$ the $(k+1)$-dimensional plane
generated by $l_{1}$ and $l_{2}$.
Assume that
$$C \ne \ggg^{n}_{k}(s_{2})\;.$$
In this case the inclusion $C \subset \ggg^{n}_{k}(s_{2})$
does not hold (indeed $C \in {\cal MC}^{n}_{k}(V)$) and the set
$C\setminus\ggg^{n}_{k}(s_{2})$
is not empty.

Consider a plane $l$ belonging to it
and define $s'=l\cap s_{2}$.
It is trivial that
$\dim s' \le k-1$
(if the inequality does not hold then
$l\in \ggg^{n}_{k}(s_{2})$).
The plane $s'$ contains the planes
$l\cap l_{1}$ and $l\cap l_{2}$
(indeed $l_{1}$ and $l_{2}$ are contained in $s_{2}$).
Each two planes from the collection
$l_{1},l_{2},l$ are adjacent
(these are elements of $C$)
and
$$\dim l\cap l_{1}=\dim l\cap l_{2}=k-1\;.$$
This implies that
$$s'=l\cap l_{1}=l\cap l_{2}\;.$$
Then
$$s'=l\cap l_{1}\cap l_{2} \subset l_{1} \cap l_{2}=s_{1}\;;$$
i.e. the planes $s_{1}$ and $s'$ are coincident
and $l\in \ggg^{n}_{k}(s_{1})$.
We obtain the inclusion
\begin{equation}
C \setminus \ggg^{n}_{k}(s_{2})\subset \ggg^{n}_{k}(s_{1})\;.
\end{equation}

Now consider a plane
$l'$ belonging to $\ggg^{n}_{k}(s_{2})\setminus
\ggg^{n}_{k}(s_{1})$.
It is contained in $s_{2}$ and
does not contain $s_{1}$.  The plane $l$
is generated by $s_{1}$
and a line which is not contained in $s_{2}$.
Therefore,
$$\dim l\cap l' <k-1$$
and these planes are not
adjacent; i.e. $l' \notin C$ and the
equality
$$C\cap(\ggg^{n}_{k}(s_{2})\setminus \ggg^{n}_{k}(s_{1}))
=\emptyset$$
is proved.
It could be rewrited in the following form
$$C\cap \ggg^{n}_{k}(s_{2})\subset \ggg^{n}_{k}(s_{1})\;.$$
Then equation (1.4.1) shows that
$C\subset \ggg^{n}_{k}(s_{1})$
and the condition
$C\in {\cal MC}^{n}_{k}(V)$
guarantes the fulfilment of the
inverse inclusion.
$\blacksquare$

\subsection{Proof of Theorem 1.4.1}

Let $f \in \CC^{n}_{k}(V)$
and $1<k<n-1$.
Then
$f$ and $f^{-1}$ transfer
any set $C \in {\cal MC}^{n}_{k}(V)$
to a set belonging to
${\cal MC}^{n}_{k}(V)$.
Proposition 1.4.2 shows that for any plane
$s\in \ggg^{n}_{k-1}(V)$ there exists a plane
$$f_{k-1}(s) \in\ggg^{n}_{i}(V),\;i=k-1 \mbox{ or } k+1$$
such that
$$f(\ggg^{n}_{k}(s))=\ggg^{n}_{k}(f_{k-1}(s))\;.$$
Define
$$V^{n}_{k-1}(f)=\{\;s \in \ggg^{n}_{k-1}(V)\;|
\;f_{k-1}(s) \in \ggg^{n}_{k-1}(V)\;\}$$
and prove the next lemma.

\begin{lemma}
If $V^{n}_{k-1}(f)\ne \emptyset$ then
$V^{n}_{k-1}(f)=\ggg^{n}_{k-1}(V)$.
\end{lemma}

{\bf Proof.}
Assume that the set $V^{n}_{k-1}(f)$ is not empty and
consiber a plane $s$ belonging to it.
We want to prove that  $V^{n}_{k-1}(f)$ contains each plane
$s'\in \ggg^{n}_{k-1}(V)$.
Lemma 1.4.2 shows that we can restrict ourself only to the case
when $s$ and $s'$ are adjacent.
In this case
the set
$\ggg^{n}_{k}(s)\cap \ggg^{n}_{k}(s')$
consists of unique plane (Lemma 1.4.1).
The mapping $f$ is a bijection
and the set
$$
f(\ggg^{n}_{k}(s)\cap \ggg^{n}_{k}(s'))=
f(\ggg^{n}_{k}(s))\cap f(\ggg^{n}_{k}(s'))=
\ggg^{n}_{k}(f_{k-1}(s))\cap \ggg^{n}_{k}(f_{k-1}(s'))\;.
$$
consists of unique plane too.
If
$f_{k-1}(s')\in \ggg^{n}_{k+1}(V)$
then this set
is not empty if and only if the plane $f_{k-1}(s)$
is contained in the plane $f_{k-1}(s')$.
However, in this case it contains
the infinite number of elements.  Therefore,
$s' \in V^{n}_{k-1}(f)$.
$\blacksquare$

Proposition 1.4.2 shows also that for any
$s\in \ggg^{n}_{k+1}(V)$
there exists a plane
$$f_{k+1}(s) \in \ggg^{n}_{i}(V),\;
i=k-1 \mbox{ or } k+1$$
such that
$$f(\ggg^{n}_{k}(s))=\ggg^{n}_{k}(f_{k+1}(s))\;.$$
Define
$$V^{n}_{k+1}(f)=\{\;s \in \ggg^{n}_{k+1}(V)\;|
\;f_{k+1}(s) \in \ggg^{n}_{k+1}(V)\;\}\;.$$
Then the following lemma
holds true.

\begin{lemma}
If $V^{n}_{k+1}(f)\ne \emptyset$ then
$V^{n}_{k+1}(f)=\ggg^{n}_{k+1}(V)$.
\end{lemma}

{\bf Proof.} The proof is similar to the proof of the previous
lemma.
$\blacksquare$

Let us define
$$i(f)=\left \{ \begin{array}{ll}
k-1\;\mbox{ if } \;V^{n}_{k-1}(f)=\ggg^{n}_{k-1}(V)\\
k+1\;\mbox{ if } \;V^{n}_{k-1}(f)=\emptyset
\end{array} \right.$$
and
$$j(f)=\left \{ \begin{array}{ll}
k+1\; \mbox{ if } \;V^{n}_{k+1}(f)=\ggg^{n}_{k+1}(V)\\
k-1\; \mbox{ if } \;V^{n}_{k+1}(f)=\emptyset\;.
\end{array} \right.$$
Lemmas 1.4.4, 1.4.5 guarantee that the mappings
$$f_{k-1}:\ggg^{n}_{k-1}(V)\to \ggg^{n}_{i(f)}(V)\;,$$
$$f_{k+1}:\ggg^{n}_{k+1}(V)\to \ggg^{n}_{j(f)}(V)$$
are well-defined.
Proposition 1.4.2
shows that they
are a bijection of $\ggg^{n}_{k-1}(V)$ onto $\ggg^{n}_{i(f)}(V)$
and a bijection of $\ggg^{n}_{k+1}(V)$ onto $\ggg^{n}_{j(f)}(V)$,
respectively. Moreover, $i(f)=k-1$ if and only if
$j(f)=k+1$.  In other words,
the next two cases
\begin{enumerate}
\item [(i)] $i(f)=k-1$ and $j(f)=k+1$
\item [(ii)] $i(f)=k+1$ and $j(f)=k-1$
\end{enumerate}
are realized.

\begin{lemma}
The mappings $f_{k-1}$ and $f_{k+1}$ preserve
the distance between planes.
\end{lemma}

{\bf Proof.}
Consider $f_{k-1}$
(for second mapping the proof is analogous).
Lemma 1.4.1 states that
planes
$s,s'\in \ggg^{n}_{k-1}(V)$
are adjacent if and only if the set
$\ggg^{n}_{k}(s)\cap \ggg^{n}_{k}(s')$
contains only one element.
The last condition holds if and only if the set
$$
\ggg^{n}_{k}(f_{k-1}(s))\cap \ggg^{n}_{k}(f_{k-1}(s'))
$$
consists of unique plane
(see the proof of Lemma 1.4.4).
Therefore, the planes $f_{k-1}(s)$ and
$f_{k-1}(s')$ are adjacent if and only if $s$ and
$s'$ are adjacent.
$\blacksquare$

An immediate verification shows that
the equality
$\alpha^{n}_{k-1}=\alpha^{n}_{k+1}$
holds if and only if $n=2k$.
In other words, in the case when $n\ne 2k$
we have $i(f)=k-1$ and $j(f)=k+1$;
the case (ii) could be realized
only for the case $n=2k$.

Denote by $\SC^{n}_{k}(V)$ the class of
transformations belonging to $\CC^{n}_{k}(V)$ and satisfying
condition (i).
Each $f\in \SC^{n}_{k}(V)$ induces the transformation
$f_{k-1}$ of $\ggg^{n}_{k-1}$
and the transformation
$f_{k+1}$ of $\ggg^{n}_{k+1}$.
Moreover,
$$f_{k-1}\in \SC^{n}_{k-1}(V)\mbox{ if } k>2$$
and
$$f_{k+1}\in \SC^{n}_{k+1}(V)\mbox{ if } k<n-2\;.$$
This implies the existence of
a family of
transformations $\{f_{i}\}^{n-1}_{i=1}$ of
$\ggg^{n}_{i}(V)$
such that
\begin{enumerate}
\item[---] $f_{k}=f$,
\item[---] for each $i=2,...,n-2$ the transformation $f_{i}$
belongs to the class $\SC^{n}_{i}(V)$ and induces $f_{k-1}$ and
$f_{k+1}$.
\end{enumerate}
Lemma 1.2.4 guarantees that $f_{1}$
induces $f_{n-1}$.  Then Theorem 1.3.2 shows that
$$f_{1}\in \LL^{n}_{1}(V) \mbox{ and }
f=L^{n}_{1\,k}(f_{1})\;.$$
We obtain the inclusion $$\SC^{n}_{k}(V)\subset\LL^{n}_{k}(V)\;.$$
The inverse inclusion is a consequence of Lemma 1.2.1 (see
Example 1.4.3).  Therefore, $$\SC^{n}_{k}(V)=\LL^{n}_{k}(V)\;.$$
It was noted above that for the case $n\ne 2k$ the classes
$\SC^{n}_{k}(V)$ and $\CC^{n}_{k}(V)$ are coincident; i.e.
Theorem 1.4.1 is proved for this case.

Let $n=2k$ and $f\in \CC^{2k}_{k}(V)$
be a transformation
satisfying  condition (ii).
Consider a transformation $g$ of
$\ggg^{2k}_{k}(V)$
defined by a bilinear form on $V$.
Lemma 1.2.2 shows
that
$$gf\in \SC^{2k}_{k}(V)=\LL^{2k}_{k}(V)\;.$$
This implies the required.
$\blacksquare$

\subsection{Proof of Theorem 1.2.1}

Denote by $g$ the transformation
of $\ggg^{n}_{m}(V)$
induced by $f$.
Lemma 1.2.3 shows that
we can restrict ourself only to the case when
$m>k$. The case $m<k$ could be reduced to it by tranponsing
of $f$ and $g$.

We want to prove that $g$ induces
a transformation of $\ggg^{n}_{m-1}(V)$.
For the case when $k=m-1$ it is trivial
(Lemma 1.2.3 states that $f$  is the required transformation).

Let $k<m-1$. Then
for a plane $s \in \ggg^{n}_{m-1}(V)$
consider planes $s_{1},s_{2}\in \ggg^{n}_{m}(V)$ such that
$s=s_{1}\cap s_{2}$ and define
$$h(s)=g(s_{1})\cap g(s_{2})\;.$$
We have
$$\ggg^{n}_{k}(s)=\ggg^{n}_{k}(s_{1})\cap\ggg^{n}_{k}(s_{2})$$
and
$$
f(\ggg^{n}_{k}(s))=
f(\ggg^{n}_{k}(s_{1}))\cap f(\ggg^{n}_{k}(s_{2}))=
\ggg^{n}_{k}(g(s_{1}))\cap \ggg^{n}_{k}(g(s_{2}))=
\ggg^{n}_{k}(h(s))$$
The trivial inclusion
$h(s)\subset g(s_{1})$
guarantees that the dimension of the plane $h(s)$
is not greater than $m-1$
(otherwise, $h(s)$ coincides with $g(s_{1})$ and the
mapping $f$ is not bijective).
However, $\dim h(s) \ge k$
(since $h(s)$ contains the plane $f(l)$
for each $l\in \ggg^{n}_{k}(s)$).
Let us prove that $\dim h(s) =m-1$.

If $k=m-2$ and
$\dim h(s) < m-1$
then $f$
maps any plane belonging to
$\ggg^{n}_{k}(s)$ to the plane
$h(s)\in\ggg^{n}_{k}(V)$;
i.e. $f$ is not a bijection.
Therefore, the inequality
$\dim h(s) < m-1$ does not hold for the
case when $k=m-2$ and we obtain the required
equality.

If $k<m-2$ then
consider
a plane $t\in \ggg^{n}_{m}(V)$
satisfying the condition
$$\dim t\cap s =m-2\;.$$
It is not difficult to see that
$$f(\ggg^{n}_{k}(t\cap s))=\ggg^{n}_{k}(g(t)\cap h(s))\;.$$
If $\dim h(s) < m-1$ then
the dimension of the plane $g(t)\cap h(s)$
is less than $m-2$ (if it fails then
the planes $g(t)\cap h(s)$ and
$h(s)$ are coincident and $f$ is not bijective).
Note that for any plane $l\in \ggg^{n}_{k}(t\cap s)$
the plane $f(l)$ is contained in the plane
$g(t)\cap h(s)$ and
$$\dim g(t)\cap h(s)\ge k\;.$$
If $k<m-3$ then consider a plane $t'\in \ggg^{n}_{m}(V)$
such that
$$\dim (t\cap s)\cap t' =m-3\;.$$
We repeat the construction until
planes $p\in \ggg^{n}_{k+1}(V)$ and $p'\in \ggg^{n}_{k}(V)$
satisfying the condition
$$f(\ggg^{n}_{k}(p))=\ggg^{n}_{k}(p')=\{p'\}$$
will be obtained.
Then the mapping $f$ is not bijective.
Therefore, the inequality $\dim
h(s) < m-1$ fails and $h(s) \in \ggg^{n}_{m-1}(V)$.

The mapping
$h$ is bijection
of $\ggg^{n}_{m-1}(V)$ into itself.
Show that it
is a transformation of
$\ggg^{n}_{m-1}(V)$.

For a plane $s' \in \ggg^{n}_{m-1}(V)$
consider planes $s'_{1},s'_{2}\in \ggg^{n}_{m}(V)$
such that
$s'=s'_{1}\cap s'_{2}$ and define
$$h'(s')=g^{-1}(s'_{1})\cap g^{-1}(s'_{2})\;.$$
The transformation $g^{-1}$ is
induced by $f^{-1}$ and the
arguments considered above guarantee that
$h'(s') \in \ggg^{n}_{m}(V)$ and
$$f^{-1}(\ggg^{n}_{k}(s'))=\ggg^{n}_{k}(h'(s'))\;.$$
However,
$$\ggg^{n}_{k}(s')=f(\ggg^{n}_{k}(h'(s')))=\ggg^{n}_{k}(hh'(s'))\;;$$
i.e. $hh'(s')=s'$ and
$h$ is a transformation of
$\ggg^{n}_{m-1}(V)$
induced by $f$.
Lemma 1.2.4 states that $h$ is induced by $g$.

It is trivial that $m-1<n-1$.
Consider the case when
$m-1\ne 1$.
For any two adjacent planes
belonging to $ \ggg^{n}_{m-1}(V)$
there exists a plane $s\in \ggg^{n}_{m}(V)$
containing them. The mappings $h$ and $h^{-1}$
transfer these planes to  adjacent planes
conteined in $\ggg^{n}_{m-1}(g(h))$
or $\ggg^{n}_{m-1}(g^{-1}(h))$, respectively.
Therefore, $h \in \CC^{n}_{m-1}(V)$.
Moreover,
$$h\in \SC^{n}_{m-1}(V)=\LL^{n}_{m-1}(V)$$
(see Subsection 1.4.3)
and
$$f=L^{n}_{m-1\,k}(h)\in \LL^{n}_{k}(V)\;.$$

In the case when $m-1=1$ ($m=2$) the similar
arguments show that
$$g\in \SC^{n}_{m}(V)=\LL^{n}_{m}(V)\;.$$
Then $f=L^{n}_{m\,k}(g)\in \LL^{n}_{k}(V)$.
$\blacksquare$

%% file: CHAPTER2.TEX
\chapter{Regular subsets of the Grassmannian manifolds}

In the chapter we introduce so-called
regular subsets of the Grassmannian manifolds.
These are discrete sets which could be considered
as some generalization of collections of
linearly independent lines.
In the next chapter they will be exploited to
define  the notion of an irregular subset.
The main result states that transformations of
the Grassmannian manifolds
preserving the class of regular sets are linear
or defined by bilinear forms.
To prove it we use the Fundamental Theorem of Projective
Geometry and the Chow Theorem.

\section{Definitions and basic properties}
\setcounter{equation}{0}
\setcounter{prop}{0}
\setcounter{lemma}{0}
\setcounter{theorem}{0}
\setcounter{exmp}{0}
\setcounter{rem}{0}

\subsection{Regular sets}

We say that a set $R \subset \ggg^{n}_{k}(V)$ is {\em
regular} if
planes belonging to it
are coordinate planes for some coordinate system
for $V$;
this coordinate system will be called
{\it associated} with the regular set $R$.
In other words, our set $R$ is regular if
there exists  a base $B$ for
$V$ such that
any plane belonging to $R$ is
generated by vectors from $B$.
This base is called
{\it associated}  with
$R$.

The class of regular
subsets of $\ggg^{n}_{k}(V)$ will be denoted by ${\cal
R}^{n}_{k}(V)$.

\begin{rem}{\rm
A coordinate system and a base associated with
a regular set are not uniquely defined.
$\blacksquare$
}\end{rem}

\begin{exmp}{\rm
Each regular subset of $\ggg^{n}_{1}(V)$ is
a colection of linearly independent lines.
$\blacksquare$
}\end{exmp}

A coordinate system for an $n$-dimensional vector space has
$$c^{n}_{k}=\frac{n!}{k!(n-k)!}$$
distinct $k$-dimensional coordinate planes.
Therefore, a regular subset of $\ggg^{n}_{k}(V)$ contains at most
$c^{n}_{k}$ elements.

We say that a regular set $R \subset \ggg^{n}_{k}(V)$ is  {\em
maximal} if any regular subset of $\ggg^{n}_{k}(V)$ containing
$R$ coincides with it.
The class of maximal regular subsets of $\ggg^{n}_{k}(V)$
is denoted by ${\cal MR}^{n}_{k}(V)$.

\begin{prop}
For any regular set $R\subset \ggg^{n}_{k}(V)$
there exists a maximal regular subset of
$\ggg^{n}_{k}(V)$
containing $R$.
The regular set $R$ is maximal if and only
if it contains $c^{n}_{k}$ elements.
\end{prop}

{\bf Proof.}
Denote by $R'$ the set of all
$k$-dimensional coordinate planes for some coordinate system
associated with $R$.  Then $R'\in {\cal MR}^{n}_{k}(V)$ and
$R\subset R'$.  Moreover, the regular set $R$ is maximal if and
only if it coincides with $R'$.
$\blacksquare$

For a regular set $R \subset \ggg^{n}_{k}(V)$
consider a coordinate system associated with it.
If our regular set is maximal then this coordinate system
is uniquely defined
(a base associated with a maximal regular
set is not uniquely defined).
Denote by $r^{n}_{k\,m}(R)$ the set of all
$m$-dimensional coordinate planes.  Then
$$r^{n}_{k\,m}(R) \in{\cal MR}^{n}_{m}(V)\;.$$
It is trivial that
$$r^{n}_{m\,m'}(r^{n}_{k\,m}(R))=r^{n}_{k\,m'}(R)
\;\;\;\forall\; R\in {\cal MR}^{n}_{k}(V)$$
and the mapping
$$r^{n}_{k\,m}:{\cal MR}^{n}_{k}(V)\to {\cal MR}^{n}_{m}(V)$$
defines an one-to-one correspondence between
${\cal MR}^{n}_{k}(V)$ and ${\cal MR}^{n}_{m}(V)$.

\begin{rem}{\rm
A family ${\cal H}$ of subsets of the set $\{1,...,n\}$
is called a {\it hypergraph}. The numbers $1,...,n$
and the elements of ${\cal H}$ are {\it verteces}
and {\it edges} (if each
edge of the hypergraph ${\cal H}$ contains
not greater than
two elements
then ${\cal H}$ is a graph).
Any $R\in {\cal R}^{n}_{k}(V)$
could be considered as some hypergraph ${\cal H}(R)$.
Fix a coordinate system associated with $R$.
Then the planes belonging to $R$
are the edges
and the  coordinate axes
are the verteces of ${\cal H}(R)$.
$\blacksquare$}
\end{rem}

\begin{exmp}{\rm
For the case $1<k<n-1$ there exists
a non-maximal regular subset of
$\ggg^{n}_{k}(V)$
having unique coordinate system associated with it.
Consider, for example, the regular set
$R \subset \ggg^{n}_{2}(V)$ with the graph represented
on the picture.
\begin{center}
\input{FIG1.PIC}
\end{center}
$\blacksquare$
}\end{exmp}

\subsection{Regular transformations of the Grassmanniam manifolds}

We say that a transformation $f$ of $\ggg^{n}_{k}(V)$
is  {\it regular} if
$f$ and $f^{-1}$
map any regular set onto
a regular set.
Each regular set is contained in some
maximal regular set. Therefore,
the transformation $f$ is regular if and only if
$f$ and $f^{-1}$
transfer  any maximal regular set to
a maximal regular set.
The class of regular transformations of $\ggg^{n}_{k}(V)$
will be denoted by $\RR^{n}_{k}(V)$.

\begin{exmp}{\rm
Let $\dim V=2$. Then any two lines belonging to $\ggg^{2}_{1}(V)$
generate a maximal regular set.
Therefore, each transformation of
$\ggg^{2}_{1}(V)$ is
regular. For the general case the similar
statement fails.
$\blacksquare$} \end{exmp}

\begin{exmp}{\rm
Show that
$$\LL^{n}_{k}(V)\subset \RR^{n}_{k}(V)\;.$$
Let $f\in \LL^{n}_{k}(V)$ and $R\in {\cal MR}^{n}_{k}(V)$.
Consider a base $B$ associated with $R$ and
a linear transformation $f'\in \LL(V)$ satisfying
the condition $f=L^{n}_{k}(f')$.
Then $f(R)$ is a set of all $k$-dimensional planes
generated by vectors from the base $f'(B)$;
i.e. $f(R)\in {\cal MR}^{n}_{k}(V)$.
Our inclusion is proved.}
$\blacksquare$
\end{exmp}

\begin{exmp}{\rm
Now consider a bijection $f$ of $\ggg^{n}_{k}(V)$ onto
$\ggg^{n}_{n-k}(V)$ defined by a bilinear form and
show that it maps
the classes
${\cal R}^{n}_{k}(V)$ and ${\cal MR}^{n}_{k}(V)$
onto the classes
${\cal R}^{n}_{n-k}(V)$ and  ${\cal MR}^{n}_{n-k}(V)$,
respectively.
This implies that for the case $n=2k$
we have the following inclusion
$$\FF^{2k}_{k}(V)\subset\RR^{2k}_{k}(V)\;.$$
Let $R \in {\cal MR}^{n}_{k}(V)$ and
$B$ be a base associated with
$R$.
The bijection
$\Phi^{n}_{k\,n-k}$
(see Subsection 1.2.2)
transfers $R$ to the set of
all $(n-k)$-dimensional planes generated by vectors
from the dual base $B^{*}$.
This implies that the conditions
$$R \in {\cal MR}^{n}_{k}(V) \mbox{ and }
\Phi^{n}_{k\,n-k}(R) \in {\cal MR}^{n}_{n-k}(V^{*})$$
are equivalent for
any set $R\subset \ggg^{n}_{k}(V)$.
The required statement is a concequence of
the following fact: the bijection
$f$ could be represented as the composition of $\Phi^{n}_{k\,n-k}$
and a linear isomorphism belonging to $\LL^{n}_{n-k}(V^{*},V)$
(Subsection 1.2.2).
$\blacksquare$}
\end{exmp}

The next theorem is the main result of the chapter.

\begin{theorem}
If $n\ne 2k$ then
$$\RR^{n}_{k}(V)=\LL^{n}_{k}(V)\;.$$
If $n=2k$ and $k>1$ then
$$\RR^{2k}_{k}(V)=\FF^{2k}_{k}(V)\;.$$
\end{theorem}

{\bf Proof Theorem 2.1.1 for the cases
$k=1,n-1$.}
If $k=1$  then
any regular subset of $\ggg^{n}_{k}(V)$
is a collection of linearly independent
lines. Therefore, in this case
Theorem 2.1.1 is a trivial consequence
of Theorem 1.3.1.

Let $f\in\RR^{n}_{n-1}(V)$.
Consider a bijection
$g$ of $\ggg^{n}_{n-1}(V)$ onto $\ggg^{n}_{1}(V)$
defined by
a bilinear form.
Then
$$f'=gfg^{-1}\in \RR^{n}_{1}(V)=\LL^{n}_{1}(V)$$
(Example 2.1.5).
Lemmas 1.2.1 and 1.2.2 shows that
$f= g^{-1}f'g$
induces a transformation of $\ggg^{n}_{m}(V)$
for each $m=1,...,n-2$.
Theorem 1.2.1 guarantees that
$f\in\LL^{n}_{n-1}(V)$.
$\blacksquare$

For the case when $1<k<n-1$ Theorem 2.1.1
is not so simple to prove.
In Section 2.3 we obtain the equality
$$\RR^{n}_{k}(V)=\CC^{n}_{k}(V)$$
showing that Theorem 2.1.1
is a consequence of the Chow Theorem.
To prove it we exploit properties of one
number characteristic of regular sets,
so-called the degree of inexactness (Section 2.2).

\subsection{Representation of $\ggg^{n}_{k}(s)$ as the
Grassmannian \\manifold}

Let $s\in \ggg^{n}_{m}(V)$.
If  $m>k$ then there exists the natural isomorphism
of $\ggg^{n}_{k}(s)$ onto the
Grassannian manifold $\ggg^{m}_{k}(\hat{V})$
(the plane $s$ could be
considered
as an $m$-dimensional vector space $\hat{V}$).

In the case $m<k$
a bijection of $\ggg^{n}_{k}(V)$
onto $\ggg^{n}_{n-k}(V)$
defined by a bilinear form $\Omega$ on $V$
maps $\ggg^{n}_{k}(s)$ onto $\ggg^{n}_{n-k}(s^{\perp}_{\Omega})$.
Let us consider $s^{\perp}_{\Omega}$ as an
$(n-m)$-dimensional vector space
${\overline V}$. Then $n-m>n-k$ and
$\ggg^{n}_{n-k}(s^{\perp}_{\Omega})$ could be
represented as the  Grassmannian manifold
$\ggg^{n-m}_{n-k}({\overline V})$ which is isomorphic to
$\ggg^{n-m}_{k-m}({\overline V})$.

We have constructed the isomorphism $f$
of $\ggg^{n}_{k}(s)$ onto
$$\ggg^{m}_{k}(\hat{V})\;
(\mbox{if }m>k)\;\mbox{ or }\;\ggg^{n-m}_{k-m}({\overline V})\;
(\mbox{if }m<k)$$
satisfying  the following condition:
for any $R \subset \ggg^{n}_{k}(s)$ the set
$f(R)$ is regular if and only if $R$ is regular.

We say that a regular set $R \subset \ggg^{n}_{k}(s)$
is {\it maximal in} $\ggg^{n}_{k}(s)$ if any regular set
$R' \subset \ggg^{n}_{k}(s)$ containing $R$
coincides with it. It is equivalent to
the existence of a maximal regular set
$\hat{R}\subset \ggg^{n}_{k}(V)$ such that
$R=\hat{R}\cap\ggg^{n}_{k}(s)$.  The regular
set $R$ is maximal in $\ggg^{n}_{k}(s)$ if and only if $f(R)$ is
a maximal regular subset of the respective Grassmannian manifold.

\begin{exmp}{\rm
Let $R \in {\cal MR}^{n}_{k}(V)$
and $s \in r^{n}_{k\,m}(R)$ (i.e. $s$ is an $m$-dimensional
coordinate plane for the coordinate system associated with $R$).
The isomorphism constructed above transfers the set
$$R(s)=\ggg^{n}_{k}(s)\cap R$$
to a maximal regular subset of
$\ggg^{m}_{k}(\hat{V})$
or $\ggg^{n-m}_{k-m}({\overline V})$.
Therefore,
$$|R(s)|=
\left \{ \begin{array}{ll} c^{m}_{k}
& \mbox{ if }\; m>k \\ c^{n-m}_{k-m} & \mbox{ if }\; m <k\;.
\end{array} \right.
$$
$\blacksquare$
}\end{exmp}

\begin{exmp}{\rm
Let $R \in {\cal MR}^{n}_{k}(V)$
and $s \in r^{n}_{k\,n-1}(R)$.
Let also $s' \in r^{n}_{k\,m}(R)$ be a
plane
contained in $s$ and $m<k$.
Then $R(s)\cap R(s')$
is the set of all planes belonging
to the maximal regular subset $R(s)$
of $\ggg^{n-1}_{k}(\hat{V})$ and containing the plane $s'$.
This implies that
$$|R(s)\cap R(s')|=c^{n-m-1}_{k-m}$$
(Example 2.1.6).
This equality will be used in Section 2.3.
$\blacksquare$
}\end{exmp}

\section{Degree of inexactness of regular sets}
\setcounter{equation}{0}
\setcounter{prop}{0}
\setcounter{lemma}{0}
\setcounter{theorem}{0}
\setcounter{rem}{0}
\setcounter{exmp}{0}
\setcounter{cor}{0}

\subsection{Definitions}

We say that a regular set $R \subset \ggg^{n}_{k}(V)$
is {\it
exact} if there exists only one maximal regular subset
of $\ggg^{n}_{k}(V)$ containing
it; i.e. a coordinate system associated with $R$ is uniquely
defined.
The class of exact regular subsets of $\ggg^{n}_{k}(V)$ is
denoted by ${\cal ER}^{n}_{k}(V)$.

\begin{exmp}{\rm
Any maximal irregular set is exact.
The inverse statement fails (see Example 2.1.2).
$\blacksquare$
}\end{exmp}

For a regular set $R \in {\cal R}^{n}_{k}(V)$
the number
$$ \deg(R)=\min\{|\hat R|-|R|\;|\;
\hat R \in {\cal ER}^{n}_{k}(V)\;,R \subset \hat R\;\}$$
will be called the {\it degree of inexactness} of
$R$.
It is trivial that the equality $\deg(R)=0$ holds if and only
if our regular set $R$ is exact.

\begin{prop}
Any regular transformation $f\in \RR^{n}_{k}(V)$
and any bijection $g$ of $\ggg^{n}_{k}(V)$ onto
$\ggg^{n}_{n-k}(V)$ defined by a bilinear form on $V$ preserve
the degree of inexactness; i.e. the equality
$$\deg(f(R))=\deg(g(R))=\deg(R)$$
holds for each regular set $R \subset \ggg^{n}_{k}(V)$.
\end{prop}

{\bf Proof.}
A regular set $R$ is contained
in some maximal regular set $R'$ if and only if
$f(R)$ is contained in the maximal regular set
$f(R')$.
This implies that $f$ preserves the class ${\cal ER}^{n}_{k}(V)$.
The  set $R$ is contained
in an exact regular set $R''$ if and only if
$f(R)$ is contained in the exact regular set
$f(R'')$. Therefore, $\deg(f(R))=\deg(R)$.
The similar arguments show that
the analogous equality holds for the mapping $g$.
$\blacksquare$

\begin{exmp}{\rm
In the case $k=1$ for any $R\in {\cal R}^{n}_{k}(V)$ we
have $\deg(R)= n-|R|$ and the regular
set $R$ is exact if and only if $R\in {\cal MR}^{n}_{k}(V)$.
Proposition 2.2.1 shows that the similar statement holds for
the case $k=n-1$.
$\blacksquare$}
\end{exmp}

The general case  is more complicated.
In the next subsection we consider regular sets containing
a lot of elements. It will be proved there that these sets have a
small degree of inexactness.

\subsection{Degree of inexactness of regular sets containing a
lot of elements}

We begin with introducing of some terms
which will be used in what follows.

Let $R'\in {\cal R}^{n}_{k}(V)$. Fix some
exact regular set
$R''$ containing $R'$
and satisfying the condition
$$\deg(R')=|R''|-|R'|\;.$$
There exists  unique maximal regular set $R\in {\cal
MR}^{n}_{k}(V)$ containing $R''$.  Let $\{x_{i}\}^{n}_{i=1}$ be
a base associated with it.  For any $i=1,...,n$ define
$l_{i}=l(x_{i})$ (recall that $l(x)$ is the line containing  a
vector $x$) and
$$R_{i}=R'\cap R(l_{i})\;.$$
Then $R_{i}$ is the  set of all planes belonging to $R'$ and
containing the line $l_{i}$.  Consider the plane
$$s_{i}=\bigcap_{l \in R_{i}}l\;.$$
It is easy to see that $R_{i}=R'\cap R(s_{i})$.
Now define
$$ n_{i}=\left \{
\begin{array}{ll} \dim s_{i}\;\;\; & \mbox{ if }R_{i}\ne
\emptyset\\ 0\;\;\;  & \mbox{ if }R_{i}=\emptyset \end{array}
\right. $$
and denote by $n(R')$ the number of all $i$ such that
$n_{i}=1$.

Note that
{\it the number $n(R')$ is independent from the choose
of an exact regular set $R''$ and the regular set
$R'$ is exact if and only if $n(R')=n$}.

\begin{exmp}{\rm
Let $R \in {\cal MR}^{n}_{k}(V)$
and
$s \in r^{n}_{k\,m}(R)$.
If $R'=R(s)$ and $m>k$ then $R'$ is a maximal
regular subset of $\ggg^{m}_{k}(\hat{V})$
(Subsection 2.1.3) and $n_{i}(R')=1$
for each $i$ such that the line $l_{i}$
is contained in $s$.
$\blacksquare$
}\end{exmp}

\begin{exmp}{\rm
For a maximal regular set $R \in {\cal MR}^{n}_{k}(V)$
consider two planes
$s \in r^{n}_{k\,n-1}(R)$ and
$s'\in r^{n}_{k\,2}(R)$
such that $s$ does not contain $s'$.
Assume that
$$R'=R(s)\cup R(s')\;.$$
The sets $R(s)$ and $R(s')$ are disjoint and
$$|R'|=|R(s)|+|R(s')|=
c^{n-1}_{k}+c^{n-2}_{k-2}$$
(Example 2.1.6).
There exists unique line $l_{i}$ transverse to $s$
(it is contained in $s'$).
Then $n_{j}=1$ if $i\ne j$
(Example 2.2.3)
and $s_{i}=s'$.
Therefore,
$n_{i}=2$ and the regular set $R'$ is not exact.
Consider
a plane $l\in R$ which contains $l_{i}$
and does not contain $s'=s_{i}$.
The intersection $l\cap s_{i}$ is the line
$l_{i}$. This implies that
the regular
set $R'\cup \{l\}$ is exact and $\deg(R')=1$.
$\blacksquare$
}\end{exmp}

For  any natural $k$ satisfying the inequality $1<k<n-1$
define
$$s^{n}_{k}=c^{n-1}_{k}+c^{n-2}_{k-2}\;.$$
Then the following statement holds true.

\begin{theorem}
If a regular set $R'\in {\cal R}^{n}_{k}(V)$ $(1<k<n-1)$
contains not less than $s^{n}_{k}$ elements then
$\deg(R') \le 1$ and the equality $\deg(R')=1$ holds if and
only if $R'$ is the set considered in Example $2.2.4$.
\end{theorem}

\begin{cor}
Each regular subset of $\ggg^{n}_{k}(V)$ $(1<k<n-1)$ containing
greater than $s^{n}_{k}$ elements is exact.
\end{cor}

\begin{rem}{\rm
We have
$$c^{n}_{k}=c^{n-1}_{k-1}+c^{n-1}_{k}=
c^{n-2}_{k-1}+c^{n-2}_{k-2}+c^{n-1}_{k}=s^{n}_{k}+c^{n-2}_{k-1}\;;$$
i.e.
$$s^{n}_{k}=c^{n}_{k}-c^{n-2}_{k-1}\;.$$
Corollary 2.2.1 states that an exact regular subset of
$\ggg^{n}_{k}(V)$ could be optained from a maximal regular set
by removing $c^{n-2}_{k-1}-1$
arbitrary elements.
$\blacksquare$}
\end{rem}

{\bf Proof of Theorem 2.2.1.}
If the regular set $R'$ is not exact then there exists a number
$i$ such that $n_{i}\ne 1$.
Consider the plane $s \in r^{n}_{k\,n-1}(R)$
transverse to $l_{i}$.
For each plane $l\in R'$ one of the following cases
$$l\in R_{i} \mbox{ or }l\in R(s)$$
is realized;
i.e. $R'$ could be represented as the union of two disjoint
sets
$$
R'=(R(s)\cap R')\cup R_{i}\;.
$$
We have
$$
|R(s)\cap R'|\le |R(s)|= c^{n-1}_{k}\;,
$$
$$
|R_{i}=R(s_{i})\cap R'|\le
|R(s_{i})|=c^{n-n_{i}}_{k-n_{i}}\le
c^{n-2}_{k-2}
$$
(in the case $n_{i}=0$ the set $R_{i}$
is empty; therefore, the
last inequality is a consequence of
the condition $n_{i}\ne 1$ and Remark 2.2.2).
Then the inequality
$$|R'|=|R(s)\cap R'|+|R_{i}|
\le c^{n-1}_{k}+c^{n-2}_{k-2} =s^{n}_{k}$$
shows that the
condition $|R'| \ge s^{n}_{k}$ holds if and only if
$$R(s)\cap R'=R(s)\;,$$
$$R_{i}=R(s_{i})$$
and
$n_{i}=2$.
$\blacksquare$

\begin{rem}{\rm
Let $k_{1}$ and $k_{2}$ be natural numbers satisfying the
condition $0<k_{1}\le k_{2}<k$.  Then an immediate verification
shows that
$$c^{n-k_{1}}_{k-k_{1}}\ge c^{n-k_{2}}_{k-k_{2}}\;.$$
This fact will be often exploited in the next
section.
$\blacksquare$
}\end{rem}

\begin{exmp}{\rm
Let $R \in {\cal MR}^{n}_{k}(V)$.
Assume that $n-k\le k<n-1$
and $R'=R(s)$, where
$s \in r^{n}_{k\,1}(R)$.
It will be proved in
the next section that for this case we
have $\deg(R')=2$.

A bijection
$f$ of  $\ggg^{n}_{k}(V)$ onto $\ggg^{n}_{n-k}(V)$ defined by a
bilinear form $\Omega$ on $V$ maps this set onto
the set of planes belonging to a maximal regular
subset $f(R)$ of $\ggg^{n}_{n-k}(V)$
and contained in the $(n-1)$-dimensional plane
$s^{\perp}_{\Omega}$.
Then Proposition 2.2.1 guarantees the fulfilment of
the equality $\deg(R')=2$ for the case when $R'=R(s)$,
$1<k\le n-k$ and $s \in r^{n}_{k\,n-1}(R)$.
$\blacksquare$ }\end{exmp}

\begin{theorem}
Let $R' \in {\cal R}^{n}_{k}(V)$ and $1<k<n-1$.
Then the following statements are fulfilled:
\begin{enumerate}
\item [{\rm(i)}]
if $n-k<k<n-1$ and $R'$ contains not less than
$c^{n-1}_{k-1}$ elements then $\deg(R') \le 2$ and the equality $\deg(R')=2$
holds if and only if there exist
$R\in {\cal MR}^{n}_{k}(V)$
and $s \in r^{n}_{k\,1}(R)$
such that $R'=R(s)$;
\item [{\rm(ii)}] if $1<k<n-k$ and $R'$ contains not
less than $c^{n-1}_{k}$ elements then $\deg(R') \le 2$ and the
equality $\deg(R')=2$ holds if and only if there exist
$R\in {\cal MR}^{n}_{k}(V)$ and
$s \in r^{n}_{k\,n-1}(R)$
such that $R'=R(s)$.
\item [{\rm(iii)}]
if $n=2k$ and $R'$ contains not less than
$c^{n-1}_{k}=c^{n-1}_{k-1}$ {\rm(see Remark 2.2.3)}
elements then $\deg(R') \le 2$ and the equality $\deg(R')=2$
holds if and only if there exist
$R\in {\cal MR}^{n}_{k}(V)$ and
$s\in r^{n}_{k\,m}(R)$ such that $m=1$ or $n-1$
and $R'=R(s)$;
\end{enumerate}
\end{theorem}

\begin{rem}{\rm
An immediate verification shows that
$c^{n-1}_{k}=c^{n-1}_{(n-k)-1}$
and
$c^{n-1}_{k}=c^{n-1}_{k-1}$ if and only if $n=2k$.
A bijection of  $\ggg^{n}_{k}(V)$ onto $\ggg^{n}_{n-k}(V)$
defined by a bilinear form on $V$ transfers a set satisfying the
conditions of statement (i) to a set satisfying the conditions of
statement (ii).  This implies that statement (ii)  is a
consequence of Proposition 2.2.1 and statement (i) (see Example
2.2.5).  $\blacksquare$ }\end{rem}

\begin{cor}
For any $R' \in {\cal R}^{n}_{k}(V)$
the following two statements hold true:
\begin{enumerate}
\item [---]
if $n-k \le k<n-1$ and $R'$ contains
greater than $c^{n-1}_{k-1}$ elements then $\deg(R') \le 1$;
\item [---]
if $1<k \le n-k$ and  $R'$ contains
greater than $c^{n-1}_{k}$ elements then $\deg(R') \le 1$.
\end{enumerate}
\end{cor}

In the next section we give the proof of
Theorem 2.2.2 and exploit it to prove
Theorem 2.1.1.

\section{Proof of Theorems 2.2.2 and 2.1.1}
\setcounter{theorem}{0}
\setcounter{lemma}{0}
\setcounter{prop}{0}
\setcounter{equation}{0}
\setcounter{defn}{0}
\setcounter{rem}{0}

\subsection{Lemmas}

This subsection is devoted to
prove a few lemmas
which will be used in Subsection 2.3.3 to prove Theorem 2.2.2.
We exploit the notation introduced in
first part of Subsection 2.2.2.
Remark 2.2.3 shows that we can
restrict
ourself only to the case when
$n-k \le k <n-1$ and $|R'|\ge c^{n-1}_{k-1}$.

\begin{lemma}
The condition $n_{i}=0$ holds for some number $i$
if and only if $n=2k$ and there exists a
plane $s\in r^{n}_{k\,n-1}(R)$ such that $R'=R(s)$.
\end{lemma}

{\bf Proof.}
If there exists a number $i$ such that $n_{i}=0$ then
the set $R_{i}$ is empty and
each plane belonging to $R'$ is contained in the
plane $s\in r^{n}_{k\,n-1}(R)$ transverse to the
line $l_{i}$; i.e. $R' \subset R(s)$. This
implies that
$$
c^{n-1}_{k-1} \le |R'| \le |R(s)|=c^{n-1}_{k}\;.
$$
An immediate verification shows that
$$c^{n-1}_{k-1} \ge c^{n-1}_{k}
\mbox{ if }n-k\le k$$
and this inequality can be replaced
by an equality if and
only if $n=2k$.
Therefore, in our case $n=2k$
and  $R'=R(s)$.
$\blacksquare$

\begin{rem}{\rm
Lemma 2.3.1 and Example 2.2.3 show that
if $n_{i}=0$ for some  $i$
then $n(R')=n-1$.}
$\blacksquare$
\end{rem}

\begin{lemma}
The inequality $n_{i} \le n-k$ holds for
any $i=1,...,n$.
\end{lemma}

\begin{rem}{\rm
For the general case the strict inequality $n_{i} <n-k$
fails.
Assume that $n=2k$ and
$R'=R(s)\cup\{l\}$,
where $s\in r^{n}_{k\,n-1}(R)$ is transverse to
the line $l_{i}$ and the plane $l\in R'$ contains $l_{i}$.
Then $R_{i}=\{l\}$ and
$n_{i}=k=n-k$.
$\blacksquare$
}\end{rem}

{\bf Proof.}
The case $n_{i}=0$ is trivial.
In the case $n_{i}>0$ the set $R_{i}$ is not empty and
there
exists a plane $l\in R'$ containing
$l_{i}$.  Then $n_{i} \le k$ and for the case $n=2k$
($n-k=k$)
the required inequality is proved.

In the case when $n-k<k$ consider the plane $s \in
r^{n}_{k\,n-1}(R)$ transverse to $l_{i}$.  Then \begin{equation}
c^{n-1}_{k-1}\le |R'| \le |R(s_{i})|+ |R(s)|=
c^{n-1}_{k}+c^{n-n_{i}}_{k-n_{i}}
\end{equation}
(see the proof of Theorem 2.2.1).
Remark 2.2.2 shows that
if $n_{i} \ge n-k+1$ then
$$c^{n-n_{i}}_{k-n_{i}}\le c^{k-1}_{2k-n-1}$$
and we can rewrite inequality (2.3.1) in the following form
\begin{equation}
c^{n-1}_{k-1}-c^{n-1}_{k}\le c^{k-1}_{2k-n-1}\;.
\end{equation}
We have
$$c^{n-1}_{k-1}-c^{n-1}_{k}=
\frac{(n-1)!}{(k-1)!(n-k)!}-\frac{(n-1)!}{k!(n-k-1)!}=
\frac{(n-1)!(2k-n)}{k!(n-k)!}=$$
$$(2k-n)\underbrace{k(k+1)...(n-2)(n-1)}_{n-k}
\frac{(k-1)!}{k!(n-k)!}\;,$$
$$c^{k-1}_{2k-n-1}= \frac{(k-1)!}{(2k-n-1)!(n-k)!}=$$
$$(2k-n)\underbrace{(2k-n+1)...(k-1)k}_{n-k}
\frac{(k-1)!}{k!(n-k)!}\;.$$
The condition $n-k < k<n-1$ guarantees that
$$(2k-n)\underbrace{k(k+1)...(n-2)(n-1)}_{n-k} >
(2k-n)\underbrace{(2k-n+1)...(k-1)k}_{n-k}$$
and inequality (2.3.2) does not hold. Therefore,
$n_{i}\le n-k$.
$\blacksquare$

\begin{lemma}
If there exists a number $i$ satisfying the
condition $n_{i} \ge 3$ then $n_{j}=1$ for any
$j \ne i$.
\end{lemma}

{\bf Proof.}
Consider the planes $s\in r^{n}_{k\,n-1}(R)$
transverse to $l_{i}$
and $s'\in r^{n}_{k\,n-2}(R)$
transverse to the two-dimensional plane generated by
$l_{i}$ and $l_{j}$.
For each plane $l\in R'$ the next
three cases could be realized.
\begin{enumerate}
\item[---] $l \in R_{i}$.
\item[---] If $l \notin R_{i}$ and $l \in R_{j}$ then the
plane $l$ does not contain the line $l_{i}$ and $l \in R(s)$.
Then the condition $l \in R_{j}$ shows that
$l \in R(s)\cap R(s_{j})$.
\item[---] If $l \notin R_{i}\cup R_{j}$ then
$l$ does not contain the lines $l_{i}$, $l_{j}$
and we have $l \in R(s')$.
\end{enumerate}
Then
$$ R'\subset R(s_{i})\cup (R(s)\cap R(s_{j}))\cup R(s')$$
and we obtain the following inequality
\begin{equation}
|R'| \le |R(s_{i})|+ |R(s)\cap R(s_{j})|+|R(s')|\;.
\end{equation}
The set
$R(s)\cap R(s_{j})$
is not empty if and only if
the plane $s_{j}$ is contained in $s$.
In this case
$$|R(s)\cap R(s_{j})|=c^{n-n_{j}-1}_{k-n_{j}}$$
(see Example 2.1.7)
and equations (2.3.3) shows that
$$
c^{n-1}_{k-1}\le c^{n-n_{i}}_{k-n_{i}}+c^{n-n_{j}-1}_{k-n_{j}}+c^{n-2}_{k}
$$
(Example 2.1.6).
For the case when $n_{i} \ge 3$ and $n_{j} \ge 2$ we have
$$c^{n-n_{i}}_{k-n_{i}}\le c^{n-3}_{k-3}$$
and
$$c^{n-n_{j}-1}_{k-n_{j}}\le c^{n-3}_{k-2}$$
(Remark 2.2.2); i.e.
$$
c^{n-1}_{k-1}  \le c^{n-3}_{k-3}+ c^{n-3}_{k-2} + c^{n-2}_{k}\;.
$$
The equality
$$c^{n-1}_{k-1}=c^{n-2}_{k-1}+c^{n-2}_{k-2}=c^{n-2}_{k-1}
+c^{n-3}_{k-3}+ c^{n-3}_{k-2}$$
implies that the last inequality could be rewrited
in the next form
$$c^{n-2}_{k-1} \le c^{n-2}_{k}\;.$$
An immediate verification shows that
it does not hold for the case
$n-k \le k$.

Then the condition $n_{i} \ge 3$ guarantees
the fulfilment of the inequality $n_{j}\le 1$.
Remark 2.3.1 implies that $n_{j}>0$ and
we obtain the required equality.
$\blacksquare$

\begin{lemma}
If the condition $n_{i}=2$ holds for some
number $i$ then there exists unique line $l_{j}$ {\rm(}$j\ne
i${\rm)} contained in the plane $s_{i}$ and such that
$n_{j}=1$.
\end{lemma}

{\bf Proof.} The equality $n_{i}=2$
implies the
existence of unique
$l_{j}$ ($j\ne i$) contained
in $s_{i}$.
The trivial inclusion $R_{i} \subset R_{j}$
shows that $s_{j}\subset s_{i}$
and $0<n_{j}\le 2$.
Therefore, if $n_{j}\ne 1$
then the planes $s_{i}$ and $s_{j}$ are
coincident.

Consider the plane
$s \in r^{n}_{k\,n-2}(R)$ transverse to $s_{i}$.
If $l\in R'$ and
$l \notin R_{i}$ then $l$
does not contain the lines $l_{i}$ and $l_{j}$;
i.e. $l \in R(s)$. In other words, for
each $l\in R'$ one of the next two cases
$$l \in R_{i}=R_{j} \mbox{ or } l \in R(s)$$
is realized.
This implies the inclusion
$R' \subset R(s)\cup R_{i}$
showing that
$$
c^{n-1}_{k-1}=c^{n-2}_{k-1}+c^{n-2}_{k-2}\le |R'| \le |R(s)|+ |R(s_{i})|=
c^{n-2}_{k}+c^{n-2}_{k-2}\;;
$$
i.e.
$$c^{n-2}_{k-1}\le c^{n-2}_{k}\;.$$
It was noted above (see the proof of Lemma 2.3.3) that for the
case $n-k \le k$ this inequality does not hold.
Therefore,
the equality $n_{j}=2$ fails and
we have $n_{j}=1$.
$\blacksquare$

Remark 2.3.1 and Lemma 2.3.3 show that in the case
when $n_{i}\ne 1,2$ for some number $i$ we
have $n(R')=n-1$.  Let us consider the case when $0< n_{i} \le 2$
for each $i=1,...,n$.

\begin{lemma}
If $0< n_{i} \le 2$ for any $i=1,...,n$
then
$$n(R')>k \mbox{ or } n(R')=1\;.$$
Moreover, the equality $n(R')=1$
holds if and only if there exists
a number $j$ such that $R'=R(l_{j})$.
\end{lemma}

\begin{rem}{\rm
Unlike to all previous lemmas considered in
Section 2.3.1, Lemma 2.3.5 will be proved for
the case $n-k\le k\le n-1$.
$\blacksquare$
}\end{rem}

\subsection{Proof of Lemma 2.3.5.}

First of all note that if
$R'=R(l_{j})$ then $n_{j}=1$
and for any $i\ne j$
the plane $s_{i}$ is generated by $l_{i}$ and $l_{j}$;
i.e. $n_{i}=2$. Therefore, for this case
we have $n(R')=1$.

Now consider the case when $k=n-1$.
Then $c^{n-1}_{k-1}=n-1$ and $R'$ is a regular
subset of $\ggg^{n}_{n-1}(V)$ containing not less than $n-1$
elements.  If $|R'|=n$ then it is a maximal
regular set and $n(R')=n$.  In the case $|R'|=n-1$
there exists a plane $l\in R$ such that $R'=R\setminus
\{l\}$.
Consider unique line $l_{j}$ transferse to $l$,
it is trivial that $R'=R(l_{j})$.

Let $n-k\le k <n-1$. Then fix a number $i$ such that $n_{i}=2$
and consider the plane $s\in r^{n}_{k\,n-1}(R)$
transverse to the line $l_{i}$. Recall that the sets
$R(s)$ and $R(s_{i})$ are disjoint
and
$$|R'|=|R'\cap R(s)|+|R_{i}|$$
(see the proof of Theorem 2.2.1).
Therefore,
$$|R'\cap R(s)| =|R'| - |R_{i}|\ge |R'| - |R(s_{i})|\ge
c^{n-1}_{k-1}-c^{n-2}_{k-2}=c^{n-2}_{k-1}$$
and
$${\hat R}=R'\cap R(s)$$
is
a regular subset of $\ggg^{n-1}_{k}({\hat V})$ containing
not less than $c^{(n-1)-1}_{k-1}$ elements
(here we consider the plane $s$ as an
$(n-1)$-dimensional vector space ${\hat V}$).

\begin{lemma}
If $\hat R$ contains $c^{n-2}_{k-1}$ elements then
$R_{i}=R(s_{i})$.
\end{lemma}

{\bf Proof.}
In this case we have
$$
|R_{i}|=|R'|-|{\hat R}|=
c^{n-1}_{k-1}-c^{n-2}_{k-1}=c^{n-2}_{k-2}\;. $$
Then the
required statement is a consequence of
the equality $|R(s_{i})|=c^{n-2}_{k-2}$
and
the trivial inclusion $R_{i} \subset R(s_{i})$.
$\blacksquare$

For the case when $n-k=1$ Lemma 2.3.5
was proved above.
Assume that it holds for the case when
$n-k<m$
and consider the case $n-k=m$.
The inductive hypothesis and the remark made before Lemma
2.3.5 guarantee
that
$$n({\hat R})>k \mbox{ or } n({\hat R})=1\;.$$
The inequality
$n(R')\ge n({\hat R})$
implies the fulfilment of the required statement
for first case.

In second case there exists a number
$j_{1}\ne i$ such that
\begin{equation}
{\hat R}=R(s)\cap R(l_{j_{1}})\;.
\end{equation}
Therefore,
$|{\hat R}|=c^{n-2}_{k-1}$
(Example 2.1.7)
and Lemma 2.3.6
implies that $R_{i}=R(s_{i})$.
Lemma 2.3.4 guarantees the existence of
unique line
$l_{j_{2}}$ ($j_{2}\ne i$) contained in $s_{i}$
and such that
$n_{j_{2}}=1$.
An immediate verification shows that
$$R(l_{j_{2}})=R(s_{i})\cup (R(s)\cap R(l_{j_{2}}))\;.$$
Then equation (2.3.4) and the equality
$R(s_{i})=R_{i}$ imply that
$$R(l_{j_{2}})=R_{i}\cup \hat{R}=R'$$
if $j_{1}=j_{2}$.

Let us consider the case when
$j_{1}\ne j_{2}$ and show that in this case
$n_{j}=1$ for each $j\ne i$;
i.e. $n(R')=n-1$.
For $j=j_{1}$ or $j_{2}$ it is trivial.
If $j\ne j_{1},j_{2}$ then denote by  $\hat{s}_{j}$
the intersection
of all planes belonging to $\hat{R}$ and
containing $l_{j}$. It is the two-dimensional plane
generated by the lines $l_{j}$ and $l_{j_{1}}$.
The condition $j_{1}\ne j_{2}$ implies the existence of
a plane
$l\in R(s_{i})\subset R'$
which contains $l_{j}$ and does not contain $l_{j_{1}}$.
The intersection of $l$
with $\hat{s}_{j}$ is the line $l_{j}$.
Therefore,
$n_{j}=1$.
$\blacksquare$

\subsection{Proof of Theorem 2.2.2}

It was noted above that we can restrict ourself
only to the case when $n-k \le k <n-1$.
The results obtained in Subsections 2.3.1 and 2.3.2
imply that
we have to consider the following four cases.
\begin{enumerate}
\item[(i)] The inequality $n_{i}>0$ holds for any $i=1,...,n$
and there exists a number $j$ such that $n_{j} \ge 3$. Then
$n(R')=n-1$ (Lemma 2.3.3).
\item[(ii)] $0< n_{i} \le 2$ for each $i=1,...,n$ and
$n(R')>k$.
\item[(iii)] $0< n_{i} \le 2$ for each $i=1,...,n$ and
$n(R')=1$. In this case Lemma 2.3.5 implies the existence of a
number $j$ such that $R'=R(l_{j})$.
\item[(iv)] The condition $n_{i}=0$ holds for some
number $i$. Then $n=2k$ and
there exists $s\in r^{n}_{k\,n-1}(R)$ such
that $R'=R(s)$ (Lemma 2.3.1).
\end{enumerate}

{\it Case }(i).
Lemma 2.3.2 implies that
$n_{j} \le n-k$
and there exist $k-1$ numbers $i_{1},...,i_{k-1}$ such that
the plane $s_{j}$ does not contain the lines
$l_{i_{1}},...,l_{i_{k-1}}$. Denote by  $l$ the
plane
generated by the lines
$l_{i_{1}},...,l_{i_{k-1}}$ and $l_{j}$.
For $i\ne j$ we have $n_{i}=1$ and
the intersection $l \cap s_{j}$ is the line $l_{j}$.
Therefore,
the regular set $R' \cup \{l\}$ is exact
and $\deg(R')=1$.

{\it Case }(ii).
Consider all numbers
$i_{1},...,i_{m}$ such that
$$n_{i_{1}}=...=n_{i_{m}}=2\;.$$
Then $n(R')=n-m$ and the condition $n(R')>k$ shows that
$m <n-k \le k$.
Denote by $s$ the plane generated by
the planes $s_{i_{1}},...,s_{i_{m}}$.
It is easy to see that $\dim s \le 2m$
and
$$n-2m =(n-m)-m> k-m> 0\;.$$
The last inequality guarantees the existence of
$k-m$ numbers $j_{1},...,j_{k-m}$ such that
$$n_{j_{1}}=...=n_{j_{k-m}}=1$$
and $s$  does not contain  the lines
$l_{j_{1}},...,l_{j_{k-m}}$.
Denote by $l$ the
plane
generated by the lines
$$l_{i_{1}},...,l_{i_{m}}, l_{j_{1}},...,l_{j_{k-m}}\;.$$
For any number $j=1,...,m$ the plane
$s_{i_{j}}$ is generated by two lines, one of them coincides with
$l_{i_{j}}$, other line $l_{i(j)}$ satisfies the conditions
$n_{i(j)}=1$ (Lemma 2.3.4) and
$$i(j)\ne j_{1},...,j_{k-m}\;.$$
Therefore, the intersection of
$l$ with each $s_{i_{j}}$ is the line $l_{i_{j}}$
and the regular set  $R' \cup \{l\}$ is exact;
i.e. $\deg(R')=1$.

{\it Case }(iii).
In this case there are $n-1$ numbers $i$
satisfying the condition $n_{i}=2$. The inequality $k \le
n-2$ shows that for a plane
$l$ belonging to the set $R\setminus R'=R(s)$
(here $s\in r^{n}_{k\,n-1}(R)$ is transverse to
the line $l_{j}$)
there exists
a number $i$ such that $n_{i}=2$ and $l$ does not contain
$l_{i}$. This implies that the regular set $R' \cup \{l\}$
is not exact and $\deg(R')\ge 2$.

Let us prove the inverse inequality. Fix $k$ numbers
$i_{1},...,i_{k}$ such that
$$n_{i_{1}}=...=n_{i_{k}}=2$$
and denote by $l$ the  plane generated by
the lines $l_{i_{1}},...,l_{i_{k}}$.
For for each
$p=1,...,k$ the intersection
$l \cap s_{i_{p}}$ is the line $l_{i_{p}}$,
indeed the plane $s_{i_{p}}$
is generated by   $l_{j}$
and $l_{i_{p}}$.
Therefore,
$$n(R' \cup \{l\})=k+1$$
and
the regular set $R' \cup \{l\}$ satisfies the conditions
of case (ii). Then there exists $l'\in R$
such that the regular set
$R' \cup \{l,l'\}$ is exact and
$\deg(R')=2$.

{\it Case }(iv).
Let $f$ be a transformation of $\ggg^{2k}_{k}(V)$
defined by a bilinear form on $V$.
Lemma 1.2.2 shows that the regular set $f(R')$
satisfies the conditions of case (iii) and
the equality $\deg(R')=2$ is a consequence of
Proposition 2.2.1.
$\blacksquare$

\subsection{Proof of Theorem 2.1.1}

Now we exploit Theorem 2.2.2 to prove the inclusion
$$\RR^{n}_{k}(V)\subset \CC^{n}_{k}(V)$$
for the case when $1<k<n-1$. The inverse inclusion
is a trivial consequence of the Chow Theorem
(see Examples 2.1.4 and 2.1.5).

Let $f\in \RR^{n}_{k}(V)$ and
$l,l' \in \ggg^{n}_{k}(V)$ be
adjacent planes. Show that the planes
$f(l)$ and $f(l')$ are adjacent.
Consider a maximal regular set $R\subset \ggg^{n}_{k}(V)$
containing $l$ and $l'$.  Then the maximal regular set
$f(R)$
contains $f(l)$ and $f(l')$;
in what follows this set will be denoted by
$R_{f}$.
For planes
$s\in r^{n}_{k\,m}(R)$
and $s'\in r^{n}_{k\,m}(R_{f})$
define
$$R(s)=R\cap \ggg^{n}_{k}(s)\;,$$
$$R_{f}(s')=R_{f}\cap \ggg^{n}_{k}(s')\;.$$
Then the next lemma holds true.

\begin{lemma}
For each plane $s \in r^{n}_{k\,n-1}(R)$ the following
statements are fulfilled:
\begin{enumerate}
\item [{\rm (i)}] if $n \ne 2k$ then there exists
a plane $s'\in r^{n}_{k\,n-1}(R_{f})$
such that
\begin{equation}
f(R(s))=R_{f}(s')\;;
\end{equation}
\item
[{\rm (ii)}] if $n=2k$ then there exists
a plane
$s'\in r^{n}_{k\,m}(R_{f})$ such that
$m=1$ or $n-1$ and equality
{\rm(2.3.5)} holds.
\end{enumerate}
\end{lemma}

{\bf Proof.}
If $k\le n-k$ then
$$\deg(R(s))=\deg(f(R(s)))=2$$
and our statement is a trivial
consequence of Theorem 2.2.2.  In the case
$k>n-k$ consider the  line $l \in  r^{n}_{k\,1}(R)$ transverse to
the plane $s$.
Then
$$\deg(R(l))=\deg(f(R(l)))=2$$
and Theorem 2.2.2 guarantees
the existence of $l'\in r^{n}_{k\,1}(R_{f})$  such that
$$f(R(l))=R_{f}(l')\;.$$
For the plane $s'\in r^{n}_{k\,n-1}(R_{f})$
transverse to the line $l'$
we have
$$f(R(s))=f(R\setminus R(l))=R_{f}\setminus R_{f}(l')=
R_{f}(s')\;.$$
Lemma 2.3.7 is proved.
$\blacksquare$

\begin{lemma}
For each plane $s \in r^{n}_{k\,k+1}(R)$ the following
statements hold true:
\begin{enumerate}
\item [{\rm (i)}] if $n \ne 2k$ then there exists
a plane $s'\in r^{n}_{k\,k+1}(R_{f})$
such that equality {\rm(2.3.5)} holds;
\item[{\rm (ii)}] if $n=2k$ then there exists
a plane
$s'\in r^{n}_{k\,m}(R_{f})$ such that
$m=k-1$ or $k+1$ and equality
{\rm(2.3.5)} holds.
\end{enumerate}
\end{lemma}

{\bf Proof.}
Let
$r^{n}_{k\,1}(R)=\{l_{i}\}^{n}_{i=1}$
and
$r^{n}_{k\,1}(R_{f})=\{l'_{i}\}^{n}_{i=1}$.
Denote by $s_{i}$ and $s'_{i}$ the
planes belonging to the sets
$r^{n}_{k\,n-1}(R)$, $r^{n}_{k\,n-1}(R_{f})$
and transverse to the lines $l_{i}$ and
$l'_{i}$, respectively.
Lemma 2.3.7 states that in the case $n\ne 2k$ for any
$i=1,...,n$ there exists a number $j_{i}$ such that
\begin{equation}
f(R(s_{i}))=R_{f}(s'_{j_{i}})\;.
\end{equation}
Let us prove statement (i) for the case when  $s$ is generated
by the lines $l_{1},...,l_{k+1}$
(for other planes belonging to
the set $r^{n}_{k\,k+1}(R)$ the proof is similar).
We have $s=\cap^{n}_{i=k+2}s_{i}$ and
$$R(s)=\bigcap^{n}_{i=k+2}R(s_{i})\;.$$
The last equality and equation (2.3.6) show that
$$f(R(s))=\bigcap^{n}_{i=k+2}R_{f}(s'_{j_{i}})=R_{f}(s')\;;$$
where $s'$ is the plane generated by the lines
$l'_{j_{1}},...,l'_{j_{k+1}}$.

Lemma 2.3.7 guarantees that
for the case $n=2k$
the following two cases could be realized:
\begin{enumerate}
\item [(a)] there exists a number $j_{1}$
such that equation (2.3.6)
holds for $i=1$;
\item [(b)] there
exists a number $j$ such that
$f(R(s_{1}))=R_{f}(l'_{j})$.
\end{enumerate}
Prove that in  case (a)
for any $i=2,...,n$ there
exists a number $j_{i}$ such that  equation (2.3.6) holds.
Then the proof of statement (ii) will be similar
to the proof of statement (i).

Assume that there exist numbers $i$ and $j_{i}$ such that
$$f(R(s_{i}))=R_{f}(l'_{j_{i}})$$
and consider the plane $\hat s = s_{1} \cap s_{i}$.
Then
$$R(\hat s)=R(s_{1}) \cap R(s_{i})$$
and
$$f(R(\hat s))=R_{f}(s'_{j_{1}})\cap R_{f}(l'_{j_{i}})\;.$$
We have
$$|f(R(\hat s))|=c^{2k-2}_{k-1}$$
(Example 2.1.7) and
$$|R(\hat s)|=c^{2k-2}_{k}$$
(Example 2.1.6).
However,
$$c^{2k-2}_{k} \ne c^{2k-2}_{k-1}$$
(see the proof of Lemma 2.3.3).
Our hypothesis fails and statement (ii)
is proved for case (a).

For case (b)
consider the transformation
$g=F^{2k}_{k\,k}(\Omega)$ of
$\ggg^{2k}_{k}(V)$
defined by some bilinear
form $\Omega$ on $V$.
Lemma 1.2.2 shows that the regular
transformation $gf$ satisfies the conditions of case (a);
i.e. it maps $R(s_{1})$ onto the set of planes
belonging to  the maximal regular set $gf(R)$
and contained in some plane
$$\overline{s}\in r^{2k}_{2k-1}(gf(R))\;.$$
Then there exists a plane
$$s''\in r^{2k}_{k\,k+1}(gf(R))$$
such that
$$gf(R(s))= gf(R)\cap\ggg^{2k}_{k}(s'')\;.$$
Lemma 1.2.2 guarantees that the plane
$$s'=(F^{2k}_{k-1\,k+1}(\Omega))^{-1}(s'') \in
r^{2k}_{k\,k-1}(R_{f})$$
satisfies the required condition.
$\blacksquare$

Now we can prove Theorem 2.1.1.
Since the planes $l$ and $l'$
are ajacent, there exists a plane $s \in r^{n}_{k\,k+1}(R)$
containing $l$ and $l'$.
Lemma 2.3.8 implies that the planes
$f(l)$ and $f(l')$
are ajacent.
It is trivial that the inverse
transformation $f^{-1}$ is regular
and the similar arguments show that
the planes
$f^{-1}(l)$ and $f^{-1}(l')$
are ajacent too.
$\blacksquare$

%% file: CHAPTER3.TEX
\chapter{Irregular subsets of the Grassmannian
manifolds}

First of all we give the definition of
irregular and maximal irregular subsets of
the Grassmannian manifolds
(Section 3.1). It will be based
on the notion of regular sets.
In Section 3.2 we introduce two number characteristics
of subsets of the Grassmannian manifolds and
use them to study of
structural properties of irregular sets.
The main result of the chapter states that
there exist maximal irregular sets which are not similar
(we say that two subsets of the Grassmannian manifold
are similar if there exists a regular transformation
maps one of these sets onto other set).

\section{Definition, examples and elementary \\properties}
\setcounter{equation}{0}
\setcounter{prop}{0}
\setcounter{lemma}{0}
\setcounter{theorem}{0}
\setcounter{exmp}{0}
\setcounter{rem}{0}

\subsection{Irregular and maximal irregular subsets}

A set $I \subset \ggg^{n}_{k}(V)$ is called {\em irregular}
if it is not regular and does not contain maximal regular
subsets.
The class of irregular subsets of $\ggg^{n}_{k}(V)$
will be denoted by ${\cal I}^{n}_{k}(V)$.

It is trivial that
if $J\subset I \in {\cal I}^{n}_{k}(V)$
and $J\notin {\cal R}^{n}_{k}(V)$ then
$J \in {\cal I}^{n}_{k}(V)$.

\begin{exmp}{\rm
For any plane $s\in \ggg^{n}_{m}(V)$
the set $\ggg^{n}_{k}(s)$ is irregular.
$\blacksquare$
}\end{exmp}

We say that an irregular set $I \subset \ggg^{n}_{k}(V)$
is {\em maximal}  if any irregular subset
of $\ggg^{n}_{k}(V)$ containing $I$ coincides with it.
The class of maximal irregular subsets of $\ggg^{n}_{k}(V)$
is denoted by ${\cal MI}^{n}_{k}(V)$.

\begin{prop}
For any irregular set $I \subset\ggg^{n}_{k}(V)$
there exists a maximal irregular subset
of $\ggg^{n}_{k}(V)$
containing it.
\end{prop}

{\bf Proof.}
Denote by ${\cal I}$ the family of all irregular subsets of
$\ggg^{n}_{k}(V)$ containing $I$ and
show that for each
linearly ordered family ${\cal J}\subset {\cal I}$
the set
$$U_{{\cal J}}=\bigcup_{U \in {\cal J}}U $$
is irregular; in other words, the family  ${\cal J}$
is bounded above.
Then Proposition 3.1.1 is a consequence of the Zorn Lemma
(see Remark 1.4.2).

Assume that the set $U_{{\cal J}}$ is not irregular
and contains a maximal regular set $R$.
For any plane $l \in R$ there exists a set
$U(l) \in {\cal J}$ containing $l$.
The family ${\cal J}$
is linearly ordered, this implies the existence of a plane $l'\in
R$ such that
$$U(l)\subset U(l')\;\;\;\forall\;l\in R\;.$$
Then $U(l')$
contains the maximal regular set $R$ and is not irregular;
i.e. our hypothesis fails and $U_{J}\in {\cal I}^{n}_{k}(V)$.
$\blacksquare$

The next simple lemmas will be  often exploited in what follows.

\begin{lemma}
For any bijection $f$
of $\ggg^{n}_{k}(V)$ onto $\ggg^{n}_{m}(V')$
$($ where $m=k$ or $n-k$$)$
the following three
conditions are equivalent:
\begin{enumerate}
\item[---] $f$ and $f^{-1}$ map
any regular  set onto a regular set;
\item[---] $f$ and $f^{-1}$ map
any irregular  set onto an irregular set;
\item[---] $f$ and $f^{-1}$ map
any maximal irregular  set onto a maximal irregular set.
\end{enumerate}
\end{lemma}

{\bf Proof.}
It is trivial.
$\blacksquare$

\begin{lemma}
An irregular set $I \subset \ggg^{n}_{k}(V)$ is maximal if and
only if for any plane $l$ belonging to  $\ggg^{n}_{k}(V) \setminus
I$ there exists a set $R \subset I$ such that
$R \cup \{l\}$ is a maximal regular subset
of $\ggg^{n}_{k}(V)$.
\end{lemma}

{\bf Proof.}
The irregular set $I$ is maximal if and only if
for any plane $l \in \ggg^{n}_{k}(V) \setminus I$ the set
$I \cup \{l\}$ is not irregular.
$\blacksquare$

\begin{exmp}{\rm
Let $s\in \ggg^{n}_{n-1}(V)$.
Then Lemma 3.1.2 implies that the irregular
set $\ggg^{n}_{1}(s)$ is maximal.

Now assume that $t\in \ggg^{n}_{1}(V)$.
In this case the bijection
of $\ggg^{n}_{n-1}(V)$ onto $\ggg^{n}_{1}(V)$
defined by some bilinear form $\Omega$ on $V$
maps $\ggg^{n}_{n-1}(t)$ onto the
maximal irregular set $\ggg^{n}_{1}(t^{\perp}_{\Omega})$
and Lemma 3.1.1 shows that the irregular set
$\ggg^{n}_{n-1}(t)$ is maximal.
$\blacksquare$
}\end{exmp}

Show that for the cases $k=1$ or $n-1$ the class
${\cal MI}^{n}_{k}(V)$ consists of  sets
considered in the last example; i.e. the following
statement holds true.

\begin{prop}
If $I\in {\cal MI}^{n}_{k}(V)$
and $k=1,n-1$ then there exists a plane
$s \in \ggg^{n}_{n-k}(V)$ such that
$I=\ggg^{n}_{k}(s)$.
\end{prop}

{\bf Proof.}
We begin with the case $k=1$.
Let $l\in \ggg^{n}_{1}(V)\setminus I$.
Then Lemma 3.1.2 implies the existence
of a  collection of lines
$l_{1},...,l_{n-1}$ belonging to $I$
and such that
$$\{l\}\cup\{l_{i}\}^{n-1}_{i=1} \in {\cal MR}^{n}_{1}(V)\;.$$
Denote by $s$  the $(n-1)$-dimensional plane generated by
$l_{1},...,l_{n-1}$.
The set $I$
does not contain
lines  transverse to $s$
(othervise, it contains a maximal  regular set).
Therefore, $I \subset \ggg^{n}_{1}(s)$.
The irregular set $I$ is maximal and
the inverse inclusion holds true.

The case $k=n-1$
could be redused to the previous case by
considering of the maximal irregular set $f(I)$,
where $f$ is
some bijection of $\ggg^{n}_{n-1}(V)$ onto $\ggg^{n}_{1}(V)$
defined by a bilinear form on $V$ (see Example 3.1.2).
$\blacksquare$

\subsection{Similar subsets of the Grassmannian manifolds}

We say that two sets $I,J\subset \ggg^{n}_{k}(V)$
are {\it similar} if there exists a regular transformation
$f\in \RR^{n}_{k}(V)$ such that $f(I)=J$.

Proposition 3.1.2 states that for the cases $k=1,n-1$
any two maximal irregular subsets of $\ggg^{n}_{k}(V)$
are similar. For the general case it fails.
In the next section we show  that
{\it
for any natural numbers $k$ and $n$ satisfying the
conditions $1<k<n-1$, $n>3$ there exist maximal irregular
subsets of $\ggg^{n}_{k}(V)$ which are not similar.}

It must be pointed out that the next simple lemma
is fulfilled.

\begin{lemma}
Let $f$ be a bijection of $\ggg^{n}_{k}(V)$ onto
$\ggg^{n}_{n-k}(V)$ defined by some bilinear form on
$V$.  Then sets $I,J\subset \ggg^{n}_{k}(V)$
are similar if and only if the sets
$f(I)$ and $f(J)$ are similar.
\end{lemma}

{\bf Proof.}
The equality $g(I)=J$ holds for some regular
transformation $g \in \RR^{n}_{k}(V)$
if and only if the regular transfornation
$$fgf^{-1}\in  \RR^{n}_{n-k}(V)$$
maps the set $f(I)$ onto the set $f(J)$.
$\blacksquare$

\subsection{Example}

Fix some plane $s\in \ggg^{n}_{m}(V)$ and
consider the sets
$$X^{n}_{k}(s)=\bigcup_{t \in \ggg^{n}_{1}(s)}\ggg^{n}_{k}(t)=
\{\;l \in \ggg^{n}_{k}\;|\;\dim l\cap s \ge 1\;\}$$
and
$$Y^{n}_{k}(s)=\bigcup_{t \in \ggg^{n}_{n-1}(s)}\ggg^{n}_{k}(t)\;.$$
For any linear transformation $f\in \LL^{n}_{k}(V)$ and any
bijection $g$ of $\ggg^{n}_{k}(V)$ onto $\ggg^{n}_{n-k}(V)$
defined by a bilinear form $\Omega$ on $V$ the next equalities
$$f(X^{n}_{k}(s))=X^{n}_{k}(L^{n}_{k,m}(s))\;,$$
$$f(Y^{n}_{k}(s))=Y^{n}_{k}(L^{n}_{k,m}(s))\;,$$
$$g(X^{n}_{k}(s))=Y^{n}_{n-k}(s^{\perp}_{\Omega})\;,$$
$$g(Y^{n}_{k}(s))=X^{n}_{n-k}(s^{\perp}_{\Omega})$$
hold true .
This implies that for any two planes
$s_{1},s_{2}\in \ggg^{n}_{m}(V)$
the sets $X^{n}_{k}(s_{1})$ and $Y^{n}_{k}(s_{1})$
are similar to
the sets $X^{n}_{k}(s_{2})$ and $Y^{n}_{k}(s_{2})$,
respectively.

Note that
$X^{n}_{k}(s)$ coincides with $\ggg^{n}_{k}(s)$
if $m=1$ and $Y^{n}_{k}(s)$ coincides with
$\ggg^{n}_{k}(s)$ if $m=n-1$.

\begin{prop}
The following statements are fulfilled:
\begin{enumerate}
\item [{\rm(i)}] if $m> n-k$ then $\ggg^{n}_{k}(V)=X^{n}_{k}(s)$;
\item [{\rm(ii)}] if $m< n-k$ then $\ggg^{n}_{k}(V)=Y^{n}_{k}(s)$;
\item [{\rm(iii)}] if $m \le n-k$ then
$X^{n}_{k}(s)\in {\cal I}^{n}_{k}(V)$;
\item [{\rm(iv)}]  if $m \ge n-k$ then
$Y^{n}_{k}(s)\in {\cal I}^{n}_{k}(V)$;
\item [{\rm(v)}] the irregular set
$X^{n}_{k}(s)$ is maximal if and only if $m=n-k$; moreover, in
this case it coincides with $Y^{n}_{k}(s)$;
\item [{\rm(vi)}] if $n\ne 2k$ and some set $I\subset
\ggg^{n}_{k}(V)$ is similar to $X^{n}_{k}(s)$ or $Y^{n}_{k}(s)$
then there exists $s'\in \ggg^{n}_{m}(V)$ such that $I$ coincides
with $X^{n}_{k}(s')$ or $Y^{n}_{k}(s')$, respectively; if $n=2k$
then the analogous statement holds only for the case when
$m=n-k$.  \end{enumerate} \end{prop}

We begin the proof of this statement with the following lemma.

\begin{lemma}
For any coordinate system for $V$ and any $m$-dimensional
plane $s$ passing through of the origin of the coordinates there
exists an $(n-m)$-dimensional coordinate plane intersecting $s$
only in the origin of the coordinates.
\end{lemma}

{\bf Proof.}
For the case $m=1$ this statement is trivial.
Assume that it holds
for any number $m$ satisfying the condition
$m< m_{0}$  and
consider the case when $m=m_{0}$.  Let $s'$ be an
$(m-1)$-dimensional plane contained in $s$.  The inductive
hypothesis implies the existence of $(n-m+1)$-dimensional
coordinate plane $t'$ intersecting $s'$ only in the origin of the
coordinates.  The intersection $s \cap t'$ is a line passing
through of the origin of the coordinates.  There exists an
$(n-m)$-dimensional coordinate plane $t$ contained in $t'$ and
intersecting this line only in the origin of the coordinates.  It
is not difficult to see that the plane $t$ satisfies the required
condition.  $\blacksquare$

{\bf Proof of Proposition 3.1.3.}
First of all note that
statements (ii) and (iv) are consequences of
statements (i), (iii) and the equalities given after
the definition of the sets $X^{n}_{k}(s)$ and
$Y^{n}_{k}(s)$.

In the case when $m> n-k$ the condition
\begin{equation}
\dim l \cap s \ge 1
\end{equation}
holds for any plane $l\in \ggg^{n}_{k}(V)$.
This implies the fulfilment of
statement (i).

Prove statement (iii).
Let $R$ be a regular set contained in
$X^{n}_{k}(s)$.  Consider a coordinate system associated with it.
Lemma 3.1.4 guarantees the existence of some $(n-m)$-dimensional
coordinate plane $t$ intersecting $s$ only in the origin of the
coordinates.  In the case $m \le n-k$ there exist $c^{n-m}_{k}\ge
1$ distinct $k$-dimensional coordinate planes contained in $t$.
These planes intersect $s$ only in the origin of the coordinates
and $X^{n}_{k}(s)$ does not contain them.  Therefore, in this
case the regular set $R$ is not maximal.

It was proved above that for the case when $m<n-k$
any regular set $R$ contained in $X^{n}_{k}(s)$
satisfies the following condition
$$|R|\le c^{n}_{k}-c^{n-m}_{k}<c^{n}_{k}-1\;.$$
Lemma 3.1.2 shows that any maximal irregular
subset of $\ggg^{n}_{k}$ contains  regular sets
containing
$c^{n}_{k}-1$ elements.
Therefore, for the case $m<n-k$
the irregular set $X^{n}_{k}(s)$ is not maximal.

In the case $m=n-k$ consider a plane
$l'$ belonging to $\ggg^{n}_{k}(V)\setminus X^{n}_{k}(s)$.
This plane intersects $s$
only in the origin of the coordinates and there exists a
coordinate system  such that $l'$ and $s$ are
coordinate planes for it.
Denote by $R$ the set of all
$k$-dimensional coordinate planes.
The equality $m=n-k$ shows that
condition (3.1.1) holds for each
plane $l$ belonging to $R\setminus \{l'\}$.
Then
$$R\setminus \{l'\} \subset X^{n}_{k}(s)$$
and Lemma 3.1.2 implies
that
$X^{n}_{k}(s)\in {\cal MI}^{n}_{k}(V)$.

Show that in the case $m=n-k$
the sets $X^{n}_{k}(s)$ and $Y^{n}_{k}(s)$
are coincident.
Any plane $l \in X^{n}_{k}(s)$
satisfies condition (3.1.1).
The equality $m=n-k$ implies the existence of an
$(n-1)$-dimensional plane containing $l$ and $s$.
Therefore, $l \in Y^{n}_{k}(s)$.
Inversely, for any $l\in Y^{n}_{k}(s)$ there
exists an $(n-1)$-dimensional plane containing $l$ and $s$.
The equality $m=n-k$ guarantees the fulfilment
of condition (3.1.1) and $l \in X^{n}_{k}(s)$.

Now prove statement (vi).
Let $f$ be a regular transformation of $\ggg^{n}_{k}(V)$
which maps $X^{n}_{k}(s)$ (or $Y^{n}_{k}(s)$)
to $I$. If $f$ is linear  then our statement
is trivial. If $f\notin \LL^{n}_{k}(V)$
then $n=2k$ and $f$ is defined
by a bilinear form on $V$. Then it maps
$X^{n}_{k}(s)$ and $Y^{n}_{k}(s)$
onto $Y^{n}_{k}(s')$ and $X^{n}_{k}(s')$,
where $s'\in \ggg^{n}_{n-m}(V)$.
The required statement is a consequence of statement (v).
$\blacksquare$

\subsection{Sets of singular restrictions of
symplectic forms}

Let $\Omega$ be a (non-singular)
symplectic form on some $n$-dimensional
vector space
$V$. Then Proposition 1.1.2 states that the number $n$ is even.
Denote by
$S^{n}_{k}(\Omega)$
the set of all planes
$s\in \ggg^{n}_{k}(V)$ such that
the restriction of the form $\Omega$
onto $s$ is singular.

Any two (non-singular) symplectic forms
$\Omega'$ and $\Omega''$
on $V$ are similar; i. e. there exists
a linear transformation $f \in \LL(V)$
such that $f^{*}(\Omega')=\Omega''$.
Then $L^{n}_{k}(f)$ maps $S^{n}_{k}(\Omega'')$
onto $S^{n}_{k}(\Omega')$ and these sets
are similar.

On odd-dimensional vector spaces there are not
non-singular symplectic forms.
Therefore,
$S^{n}_{k}(\Omega)=\ggg^{n}_{k}(V)$
if the number $k$ is odd.
In the case when $k$ is even
$S^{n}_{k}(\Omega)\in {\cal I}^{n}_{k}(V)$;
however, this
irregular set is not maximal.
We do not give the proof of this fact here.

\section{Number characteristics of irregular sets}
\setcounter{equation}{0}
\setcounter{prop}{0}
\setcounter{lemma}{0}
\setcounter{defn}{0}
\setcounter{theorem}{0}
\setcounter{rem}{0}
\setcounter{exmp}{0}
\setcounter{cor}{0}

\subsection{Definitions}

For a set $I\subset \ggg^{n}_{k}(V)$
consider the set
$$N_{1}(I)=\{\;t\in \ggg^{n}_{1}(V)\;|\;\ggg^{n}_{k}(t) \subset
I\;\}\;.$$
If this set is not empty then the lines belonging to it
generate some plane $s_{1}(I)$.
Define
$$n_{1}(I)=
\left\{
\begin{array}{ll}
\dim s_{1}(I) & \mbox{ if }\; N_{1}(I)\ne \emptyset\\
0             & \mbox{ if }\; N_{1}(I)= \emptyset\;.
\end{array}\right.$$
Now consider the dual set
$$N_{n-1}(I)=\{\;t \in \ggg^{n}_{n-1}(V)\;|\;\ggg^{n}_{k}(t)
\subset I\;\}\;.$$
If this set is not empty then the intersection
of all planes belonging to it is a plane passing through of the
origin of the coordinates
(we suppose that the
origin of the coordinates
is a zero-dimensional plane).
This plane
will be denoted by $s_{n-1}(I)$.
Define
$$n_{n-1}(I)= \left\{
\begin{array}{ll} \dim s_{n-1}(I) & \mbox{ if }\; N_{n-1}(I)\ne
\emptyset\\
n & \mbox{ if }\; N_{n-1}(I)= \emptyset\;.
\end{array}\right.$$

\begin{lemma}
For any set $I\subset \ggg^{n}_{k}(V)$
and any bijection
$f$ of $\ggg^{n}_{k}(V)$ onto
$\ggg^{n}_{n-k}(V)$ defined by a bilinear form $\Omega$
on $V$ we
have
$$n_{1}(f(I))=n-n_{n-1}(I)\;,$$
$$n_{n-1}(f(I))=n-n_{1}(I)\;.$$
\end{lemma}

{\bf Proof.}
It is not difficult to see that
$F^{n}_{1\,n-1}(\Omega)$ and $F^{n}_{n-1\,1}(\Omega)$
transfer the sets
$N_{1}(I)$ and $N_{n-1}(I)$
to the sets
$$N_{n-1}(f(I)),\;N_{1}(f(I))\;,$$
respectively.
Therefore,
$$s_{1}(f(I))=(s_{n-1}(I))^{\perp}_{\Omega}\;,$$
$$s_{n-1}(f(I))=(s_{1}(I))^{\perp}_{\Omega}$$
and we get the required.
$\blacksquare$

\begin{prop}
If $I\in {\cal I}^{n}_{k}$ then
$n_{1}(I) \le n-k$ and $n_{n-1}(I) \ge n-k$.
\end{prop}

{\bf Proof.}
Let
$t_{1},...,t_{n_{1}(V)}$ be a collection of lines belonging
to the set $N_{1}(V)$
and generating the plane $s_{1}(I)$.
Then there exist lines
$$t_{n_{1}(V)+1},...,t_{n}\in \ggg^{n}_{1}(V)$$
such that
$$T=\{t_{i}\}^{n}_{i=1}\in {\cal MR}^{n}_{1}(V)\;.$$
In the case when
$n_{1}(I)> n-k$ each plane belonging to
the maximal regular set
$$R=r^{n}_{1\,k}(T)\in {\cal MR}^{n}_{k}(V)$$
contains at least one of the lines
$t_{1},...,t_{n_{1}(V)}$.
This implies the inclusion $R \subset V$
which can not hold, since
the set $I$ is irregular.
We have proved first inequality.

Now consider the set $f(I) \in {\cal I}^{n}_{n-k}(V)$,
where $f$ is  some bijection of $\ggg^{n}_{k}(V)$
onto $\ggg^{n}_{n-k}(V)$
defined by a bilinear form on $V$.
Then $n_{1}(f(I)) \le k$ and
second inequality is a consequence of
Lemma 3.2.1.
$\blacksquare$

\begin{exmp}{\rm
Let $s\in \ggg^{n}_{m}(V)$.
Then
$$n_{1}(X^{n}_{k}(s))=m \mbox{ if } m \le n-k$$
and
$$n_{n-1}(Y^{n}_{k}(s))=m \mbox{ if } m \ge n-k\;.$$
In the case $m <n-k$ each plane $t \in \ggg^{n}_{n-1}(V)$
contains a $k$-dimensional plane intersecting
$s$ only in the origin of the coordinates.
Therefore,
$$\ggg^{n}_{k}(t) \not\subset X^{n}_{k}(s)
\;\;\;\forall\;t \in \ggg^{n}_{n-1}(V)\;;$$
i.e.
$$n_{n-1}(X^{n}_{k}(s))=0\;.$$
Then Lemma 3.2.1 guarantees the fulfilment of
the equality
$$n_{1}(Y^{n}_{k}(s))=0\mbox{ if } m > n-k\;.$$
If $m=n-k$ then the sets
$X^{n}_{k}(s)$ and $Y^{n}_{k}(s)$
are coincident and we obtain the following
equality
$$n_{1}(X^{n}_{k}(s))=n_{n-1}(X^{n}_{k}(s))=m\;.$$
In Subsection 3.2.5 we  construct
a maximal irregular set $I\in {\cal MI}^{n}_{k}(V)$ satisfying
the condition $n(I) < n-k$.
$\blacksquare$}
\end{exmp}

\begin{exmp}{\rm
If $n\ne 2k$ then each regular transformation of
$\ggg^{n}_{k}(V)$ is linear
and the equalities
$$n_{1}(I)=n_{1}(J) \mbox{ and }n_{n-1}(I)=n_{n-1}(J)$$
hold for any two similar subsets
$I,J\subset \ggg^{n}_{k}(V)$.
For the case when $n=2k$ the analogous statement fails.
Consider, for example, the similar sets
$$I=X^{2k}_{k}(s),\;s\in \ggg^{2k}_{k-1}(V)$$
and
$$J=Y^{2k}_{k}(s^{\perp}_{\Omega})=f(X^{2k}_{k}(s))\;,$$
where $f$ is the transformation of
$\ggg^{2k}_{k}(V)$
defined by
some bilinear form $\Omega$ on $V$.
For this case we have
$n_{1}(I)=k-1$ and $n_{1}(J)=0$
(see Example 3.2.1).
$\blacksquare$
}\end{exmp}

\subsection{Number characteristics and
structural properties of irregular sets I}

The inclusion $I \subset J$ guarantees the fulfilment
of the
inequalities
$$n_{1}(I) \le n_{1}(J) \mbox{ and }
n_{n-1}(I) \ge n_{n-1}(J)\;.$$
Therefore, if a set $I \subset
\ggg^{n}_{k}(V)$ contains one of the sets $X^{n}_{k}(s)$ or
$Y^{n}_{k}(s)$ then
$$n_{1}(I) \ge \dim s \mbox{ or }
n_{n-1}(I) \le \dim s\;,$$
respectively.
Now we prove that for maximal irregular sets the inverse statement
holds true.

\begin{theorem}
For any maximal regular set $I \subset \ggg^{n}_{k}(V)$
we have the following two inclusions
$$X^{n}_{k}(s_{1}(I)) \subset I\;,$$
$$Y^{n}_{k}(s_{n-1}(I)) \subset I\;.$$
\end{theorem}

{\bf Proof.}
We begin with the consideration of first inclusion.
In the case $n_{1}(I) =0$
the set $X^{n}_{k}(s_{1}(I))$ is empty and it is
trivial.
Assume that $m=n_{1}(I) \ge 1$ and
the inclusion fails.
Then  the set $X^{n}_{k}(s_{1}(I)) \setminus I$
is not empty. Let $l$ be a plane belonging to it.
The irregular set $I$ is maximal and
Lemma 3.1.2 guarantees the existence of a set
$R \subset I$
such that
$$R\cup \{l\}\in {\cal MR}^{n}_{k}(V)\;.$$
Consider a collection of lines
$t_{1},...,t_{m}\in N_{1}(I)$
generating $s_{1}(I)$.
Lemma 3.1.4 shows that for the coordinate
system associated with the maximal regular set
$R\cup \{l\}$ there
exist $n-m$ coordinate axes which are not contained in the plane
$s_{1}(I)$.  Denote them by $t_{m+1},...,t_{n}$.  Then
$$T=\{t_{i}\}^{n}_{i=1}\in {\cal MR}^{n}_{1}(V)\;.$$
Let us prove the inclusion
$$R'=r^{n}_{1\,k}(T) \subset I$$
showing that the set $I$ is not irregular.
Then first inclusion will be proved.

Consider the coordinate system associated with
the maximal regular set $R'$. Then $R'$ could be represented as
the union
$$R'=R'_{1}\cup R'_{2}\;,$$
where $R'_{1}$ is the set of
all $k$-dimensional coordinate planes containing at least one of
the lines $t_{1},...,t_{m}$ and $R'_{2}$ is the set of all
$k$-dimensional coordinate planes generated by the lines
$$t_{i_{1}},...,t_{i_{k}}\;,\;\;\;m+1 \le i_{1}<...<i_{k}\le n$$
(recall that the set $I$ is irregular and Proposition 3.2.1 shows
that $k \le n-m$).  The set $I$ contains $R'_{1}$ (since
$t_{1},...,t_{m}$ are lines belonging to $N_{1}(V)$).  It is easy
to see that
$$R'_{2}=(R\cup\{l\})\setminus X^{n}_{k}(s_{1}(I))$$
This implies that
$l\notin R'_{2}$
(indeed $l\in X^{n}_{k}(s_{1}(I))$).
Therefore, $R'_{2} \subset R \subset I$
and the inclusion  $R'\subset I$ is proved.

Let $f$ be a bijection of $\ggg^{n}_{k}(V)$ onto
$\ggg^{n}_{n-k}(V)$ defined by a bilinear form on $V$.
Then
$f(I)\in {\cal MI}^{n}_{n-k}(V)$
and Lemma 3.2.1 shows that
the inverse mapping  $f^{-1}$
transfers the inclusion
$$X^{n}_{n-k}(s_{1}(f(I)))\subset f(I)$$
to the required inclusion.
$\blacksquare$

\begin{exmp}{\rm
For irregular subsets which are not maximal
the inclusions from Theorem 3.2.1 fail.
Fix planes
$l \in \ggg^{n}_{k}(V)$
and $s \in \ggg^{n}_{n-k}(V)$
satisfying the condition
\begin{equation}
0<\dim l\cap s <n-k
\end{equation}
and consider the set
$$I=X^{n}_{k}(s)\setminus \{l\}\;.$$
Clearly, this irregular set is not maximal.
Inequality (3.2.1) implies the existence of
$n-k$ linearly independent lines
$t_{1},...,t_{n-k}$
generating $s$ and such that
the plane $l$ does
not contain each $t_{i}$.
Then
$$\ggg^{n}_{k}(t_{i})\subset I
\;\;\;\forall\;i=1,...,n-k$$
and $s_{1}(I)=s$.
It is trivial that
first inclusion from Theorem 3.2.1 does not
hold.
$\blacksquare$}
\end{exmp}

\begin{cor}
Let $I \in {\cal I}^{n}_{k}(V)$. If
\begin{equation}
n_{1}(I)=n-k
\end{equation}
then $X^{n}_{k}(s_{1}(I))$ is unique maximal irregular set
containing $I$. If
\begin{equation}
n_{n-1}(I)=n-k
\end{equation}
then $X^{n}_{k}(s_{n-1}(I))$ is unique maximal irregular set
containing $I$.
Therefore, if the irregular set $I$ is maximal and
one of conditions $(3.2.2)$ or $(3.2.3)$
holds
then there exists $s\in \ggg^{n}_{n-k}(V)$
such that
$I=X^{n}_{k}(s)$.
\end{cor}

{\bf Proof.}
It is a trivial consequence of Theorem 3.2.1 and
statement (v) of Proposition 3.1.3.
$\blacksquare$

\subsection{Number characteristics and
structural properties of irregular sets II}

Now we want to study the intersection
of arbitrary irregular set
$I\in {\cal I}^{n}_{k}(V)$ with $\ggg^{n}_{k}(t)$.
Recall that for any $t\in  \ggg^{n}_{m}(V)$
the set $\ggg^{n}_{k}(t)$
could be considered as some Grassmannian manifold
(see Subsection 2.1.3).

\begin{theorem}
For any irregular set $I \subset \ggg^{n}_{k}(V)$ the following
two statements hold true:
\begin{enumerate}
\item[{\rm(i)}]
if $m<k$  and a plane $s \in \ggg^{n}_{m}(V)$ satisfies the
condition $s \subset s_{1}(I)$ then for any plane $t\in
\ggg^{n}_{n-m}(V)$ transverse to $s$
the set $I \cap \ggg^{n}_{k}(t)$ does not contain maximal regular subsets of
$\ggg^{n}_{k}(t)$;
\item[{\rm(ii)}] if $m>k$ and a plane $s \in
\ggg^{n}_{m}$ satisfies the condition
$s_{n-1}(I) \subset s$ then
for any plane $t\in \ggg^{n}_{n-m}(V)$ transverse to $s$ the set
$I \cap \ggg^{n}_{k}(t)$
does not contain maximal
regular subsets of $\ggg^{n}_{k}(t)$.
\end{enumerate}
\end{theorem}

{\bf Proof.}
First of all we prove statement (i) for the case
when the irregular set $I$ is maximal.
In this case we have the inclusion
$X^{n}_{k}(s)\subset I$.
Assume that our statement fails; i.e. there exists a plane
$t\in \ggg^{n}_{n-m}(V)$ transverse to $s$ and such
that $I \cap \ggg^{n}_{k}(t)$ contains a
maximal irregular subset $R$ of $\ggg^{n}_{k}(t)$.
Consider $t$ as an $(n-m)$-dimensional vector space
${\hat V}$. Then $R \in {\cal MR}^{n-m}_{k}({\hat V})$.
Let
$$r^{n-m}_{k\,1}(R)=\{c_{i}\}_{i=1}^{n-m}\;.$$
Let also
$c_{n-m+1},...,c_{n}$
be a collection of lines belonging to $N_{1}(I)$
and generating the plane $s$.
Then
$$C=\{c_{i}\}^{n}_{i=1}\in {\cal MR}^{n}_{1}(V)\;.$$
Consider the coordinate system associated with
the maximal regular set
$$R'=r^{n}_{1\,k}(C)\in{\cal MR}^{n}_{k}(V)$$
This set could be represented as the union
$$R'=R\cup R''\;,$$
where $R''$ is the set of
all $k$-dimensional coordinate planes containing at least one of
the lines $c_{n-m+1},...,c_{n}$.
It is trivial that $R'' \subset X^{n}_{k}(s)$
and we obtain the inclusion
$R'\subset I$ disproving our hypothesis.

In the case when the irregular set $I$ is not maximal
consider  some maximal irregular set $I'$
containing $I$. Then our statement
is a consequence of the trivial inclusion
$$I\cap\ggg^{n}_{k}(t)\subset I'\cap\ggg^{n}_{k}(t)\;.$$

The proof of statement (ii)
is similar to the proof of second inclusion from
Theorem 3.2.1.
Let $f$ be the bijection of $\ggg^{n}_{k}(V)$ onto
$\ggg^{n}_{n-k}(V)$ defined by some bilinear form $\Omega$ on $V$.
Then
$$f(I)\in {\cal I}^{n}_{n-k}(V)
\mbox{ and }s^{\perp}_{\Omega}\subset s_{1}(f(I))\;.$$
We have the inequality $n-m < n-k$ and
for any plane
$t\in \ggg^{n}_{n-m}(V)$ transverse to $s$ the plane
$t^{\perp}_{\Omega}$ is transverse to $s^{\perp}_{\Omega}$.
Statement (i) states that the set
$$f(I)\cap\ggg^{n}_{n-k}(t^{\perp}_{\Omega})$$
does not contain maximal regular subsets of
$\ggg^{n}_{n-k}(t^{\perp}_{\Omega})$.
The mapping $f^{-1}$ transfers
this set to $I \cap \ggg^{n}_{k}(t)$,
it also maps the class of maximal regular subset of
$\ggg^{n}_{n-k}(t^{\perp}_{\Omega})$ onto the class
of  maximal regular subset of
$\ggg^{n}_{k}(t)$.
This implies the required.
$\blacksquare$

\begin{cor}
If for an irregular set $I\subset \ggg^{n}_{k}(V)$ there
exist two planes $s_{1}\in \ggg^{n}_{m_{1}}(V)$ and
$s_{2}\in\ggg^{n}_{m_{2}}(V)$ such that
$m_{1} \le n-k$, $m_{2} \ge n-k$ and
$$X^{n}_{k}(s_{1})\subset I\;,\;
Y^{n}_{k}(s_{2})\subset I$$
then $s_{1} \subset s_{2}$.
\end{cor}

{\bf Proof.}
The inclusion is a consequence of the following fact:
if a line $p\in\ggg^{n}_{1}(V)$ is not
contained in $s_{2}$ then it is not contained in $s_{1}$.
The condition $p\not\subset s_{2}$ guarantees the existence
of a plane $t \in \ggg^{n}_{n-1}(V)$
transverse to $p$ and containing $s_{2}$.
Then $\ggg^{n}_{k}(t) \subset I$
and
Theorem 3.2.2 shows that $p\not\subset s_{1}$.
$\blacksquare$

\subsection{Irregular and maximal irregular
subsets of $\ggg^{n}_{k}(t)$}

Let $t\in \ggg^{n}_{m}(V)$ and $I\subset \ggg^{n}_{k}(t)$.
We say that $I$
is an {\it an irregular} subset of $\ggg^{n}_{k}(t)$
if it is not regular
and does not contain maximal regular subsets of
$\ggg^{n}_{k}(t)$
(the definition of maximal regular subsets of $\ggg^{n}_{k}(t)$
was given in Subsection 2.1.3).
The irregular  set $I$ is called
{\it maximal} in $\ggg^{n}_{k}(t)$ if each
irregular subset of $\ggg^{n}_{k}(t)$
containing $I$ coincides with it.

In Subsection 2.1.3 we constructed the isomorphism
of $\ggg^{n}_{k}(t)$ onto some Grassmannian manifold.
It is not difficult to see that
it maps the classes of irregular and maximal irregular
subsets of $\ggg^{n}_{k}(t)$ onto the classes of
irregular and maximal irregular
subsets of the respective
Grassmannian manifold.

\begin{prop}
For any maximal irregular subset $I$ of $\ggg^{n}_{k}(t)$
there exists a maximal irregular set
$J\subset \ggg^{n}_{k}(V)$ satisfying the condition
$$J\cap \ggg^{n}_{k}(t)=I\;.$$
\end{prop}

{\bf Proof.}
We restrict ourself only to the case when $m>k$.
The case $m<k$ could be redused to it by
considering of the maximal irregular subset $f(I)$
of $\ggg^{n}_{n-k}(t^{\perp}_{\Omega})$,
where $f$ is
the bijection of $\ggg^{n}_{k}(V)$ onto $\ggg^{n}_{n-k}(V)$
defined by  some bilinear form $\Omega$ on $V$
(see the proof of Theorem 3.2.2).

First of all show that
for a plane $s\in \ggg^{n}_{n-m}(V)$
transverse to $t$ the set
$$I'=X^{n}_{k}(s)\cup I$$
is an irregular subset
of $\ggg^{n}_{k}(V)$.
If it fails then $I'$ contains a maximal regular subset
$R$ of $\ggg^{n}_{k}(V)$.
Consider the coordinate system associated with it.
It is trivial that $s$ and $t$ are coordinate
planes for this system. Then $R\cap \ggg^{n}_{k}(t)$ is a maximal
regular subset of $\ggg^{n}_{k}(t)$ contained in $I$.

For a maximal irregular set $J\subset \ggg^{n}_{k}(V)$
containing $I'$ we have
$$I\subset J\cap \ggg^{n}_{k}(t)\;.$$
Theorem 3.2.2 states that $J\cap \ggg^{n}_{k}(t)$
is an irregular subset of $\ggg^{n}_{k}(t)$
(since $s\subset s_{1}(J)$)
and we get the required equality.
$\blacksquare$

There exists a maximal irregular set
$I \subset \ggg^{n}_{k}(V)$ satisfaing the conditions
of Theorem 3.2.2 and such that
$I\cap \ggg^{n}_{k}(t)$ is not a maximal irregular
subset of $\ggg^{n}_{k}(t)$.
In what follows we show that this set
is not similar to the maximal irregular sets
considered in Subsection 3.1.3.

\begin{theorem}
For a plane $s\in \ggg^{n}_{n-k-1}(V)$ and
a plane $t\in \ggg^{n}_{k+1}(V)$ transverse to $s$
there exists a maximal irregular set
$I \subset \ggg^{n}_{k}(V)$
containing $X^{n}_{k}(s)$ and such that
$I\cap\ggg^{n}_{k}(t)$
is not a maximal irregular subset of
$\ggg^{n}_{k}(t)$.
This set
satisfies the condition
$$n_{1}(I)=n-k-1\;.$$
\end{theorem}

We have also
the following dual statement.

\begin{theorem}
For a plane $s\in \ggg^{n}_{n-k+1}(V)$ and
a plane $t\in \ggg^{n}_{k-1}(V)$ transverse to $s$
there exists a maximal irregular set
$I \subset \ggg^{n}_{k}(V)$
containing $Y^{n}_{k}(s)$ and such that
$I\cap\ggg^{n}_{k}(t)$
is not a maximal irregular subset of
$\ggg^{n}_{k}(t)$.
This set
satisfies the condition
$$n_{n-1}(I)=n-k+1\;.$$
\end{theorem}

\subsection{Proof of Theorems 3.2.3 and 3.2.4}

{\bf Proof of Theorem 3.2.3.}
We shall devide the proof
into two steps.

{\it First step.}
Let us consider some plane $s' \in  \ggg^{n}_{k+2}(V)$
containing $t$. The intersection
of $s'$ with $s$  is a line and there exists a  plane
$t'\in \ggg^{n}_{k+1}(V)$  contained in $s'$ and
such that this line is not contained in it.
It is easy to see that
$t'$ transverse to $s$ and
the planes $t$ and
$t'$ are ajacent.

Denote by $l$ the intersection of
$t$ and $t'$. Then $l\in \ggg^{n}_{k}(V)$
(since the planes $t$ and $t'$ are ajacent).
Fix two lines
$p,p'\in \ggg^{n}_{1}(V)$ contained in the planes
$t$, $t'$ (respectively) and such that
$l$ contains $p$ and does not contain $p'$.
Define
$$I'=X^{n}_{k}(s)\cup(\ggg^{n}_{k}(t')\cap\ggg^{n}_{k}(p'))
\cup(\ggg^{n}_{k}(t)\cap\ggg^{n}_{k}(p)\setminus \{l\})\;.$$
Then
$$I'\cap \ggg^{n}_{k}(t)=\ggg^{n}_{k}(t)\cap\ggg^{n}_{k}(p)
\setminus \{l\}\;,$$
$$I'\cap \ggg^{n}_{k}(t')=\ggg^{n}_{k}(t')\cap\ggg^{n}_{k}(p')$$
are irregular subsets of $\ggg^{n}_{k}(t)$
and $\ggg^{n}_{k}(t')$, respectively.
Moreover, the last set is a maximal irregular subset
of $ \ggg^{n}_{k}(t')$
(see Examples 3.1.2).

Show that $I'\in {\cal I}^{n}_{k}(V)$.
Assume that it fails and $I'$ contains a maximal regular set
$R\in {\cal MR}^{n}_{k}(V)$. Consider the coordinate
system associated with it. Then
$t$ or $t'$ is a coordinate plane for this system.
This implies (see the proof of Proposition 3.2.2)
that one of the sets $I'\cap\ggg^{n}_{k}(t)$
or $I'\cap\ggg^{n}_{k}(t')$ contains a maximal regular
subset of $\ggg^{n}_{k}(t)$ or $\ggg^{n}_{k}(t')$,
respectively.

Now consider a maximal irregular set
$I\subset \ggg^{n}_{k}(V)$ containing $I'$.
It is trivial that
$$I'\cap\ggg^{n}_{k}(t)\subset I\cap\ggg^{n}_{k}(t)\;,$$
$$I'\cap\ggg^{n}_{k}(t')\subset I\cap\ggg^{n}_{k}(t')\;.$$
We have $s\subset s_{1}(I)$
and Theorem 3.2.2 guarantees that
$I\cap\ggg^{n}_{k}(t)$ and $I\cap\ggg^{n}_{k}(t')$
are irregular subset of
$\ggg^{n}_{k}(t)$ and $\ggg^{n}_{k}(t')$.
Recall that $I'\cap\ggg^{n}_{k}(t')$
is a maximal irregular subset of $\ggg^{n}_{k}(t')$.
Therefore,
$$I\cap\ggg^{n}_{k}(t')=I'\cap\ggg^{n}_{k}(t')=
\ggg^{n}_{k}(p')\cap\ggg^{n}_{k}(t')\;.$$
Then $l \notin I$ indeed
$l\in \ggg^{n}_{k}(t')\setminus \ggg^{n}_{k}(p')$.

Any plane $l'$
belonging to
$\ggg^{n}_{k}(t)\setminus \ggg^{n}_{k}(p)$
satisfies the condition
$$\dim l\cap l'=k-1$$
There
exist a collection of $k$ lines
$p_{1},...,p_{k}$ which
generate $l'$ and are not contained in $l$.
Denote by $R'$ the set of all $k$-dimensional
planes generated by the lines
$p,p_{1},...,p_{k}$. It is a maximal regular
subset of $\ggg^{n}_{k}(t)$.
It is easy to see that
$l'\in R'$
and
$$R'\setminus \{l'\}\subset
(\ggg^{n}_{k}(p)\setminus \{l\})\cap \ggg^{n}_{k}(t)
\subset I\cap \ggg^{n}_{k}(t)\;.$$
This implies that $l'\notin I$
(since $I\cap \ggg^{n}_{k}(t)$ is an irregular subset
of $\ggg^{n}_{k}(t)$).
We obtain the equality
$$I\cap\ggg^{n}_{k}(t)=I'\cap\ggg^{n}_{k}(t)=
(\ggg^{n}_{k}(p)\setminus \{l\})\cap \ggg^{n}_{k}(t)\;.$$
Proposition 3.1.2 shows that it
is not a maximal irregular subset
of $\ggg^{n}_{k}(t)$.

{\it Second step.}
It is trivial that $n_{1}(I)\ge n-k-1$.
Shows that
\begin{equation}
X^{n}_{k}({\hat s})\not\subset I\;\;\;
\forall\;{\hat s}\in \ggg^{n}_{n-k}(V)\;.
\end{equation}
Then the required equality
will be proved.

For a plane $t\in
\ggg^{n}_{k+1}(V)$ consider the intersection of
the set $X^{n}_{k}({\hat s})$ with
$\ggg^{n}_{k}(t)$. It is trivial that
$$\dim {\hat s}\cap t \ge 1\;\;\;\forall\;
t\in\ggg^{n}_{k+1}(V)\;.$$
If $\dim {\hat s}\cap t >1$ then
any plane $l$ belonging to
$\ggg^{n}_{k}(t)$ satisfies the inequality
$$\dim {\hat s}\cap l \ge 1\;.$$
This implies that
$l\in X^{n}_{k}(s)$ and our intersection coincides with
$\ggg^{n}_{k}(t)$.

In the case when $\dim {\hat s}\cap t =1$ consider the line
$p={\hat s}\cap t$. An immediate verification
shows that
$$X^{n}_{k}({\hat s})\cap\ggg^{n}_{k}(t)=
\ggg^{n}_{k}(t)\cap\ggg^{n}_{k}(p)\;.$$
Proposition 3.1.2 guarantees
that it is
a maximal irregular subset of
$\ggg^{n}_{k}(t)$.

These arguments show that (3.2.4)
holds and Theorem 3.2.3 is proved.
$\blacksquare$

{\bf Proof of Theorem 3.2.4.}
It is trivial that
for each bilinear form $\Omega$ on $V$
the planes
$$s^{\perp}_{\Omega}\in \ggg^{n}_{k-1}(V)
\mbox{ and }t^{\perp}_{\Omega}\in \ggg^{n}_{n-k+1}(V)$$
are transverse. Theorem 3.2.3 guarantees the existence
of the maximal irregular set
$J \subset \ggg^{n}_{n-k}(V)$
such that
$$X^{n}_{n-k}(s^{\perp}_{\Omega}) \subset J$$
and
$I\cap\ggg^{n}_{n-k}(t^{\perp}_{\Omega})$
is not a maximal irregular subset of
$\ggg^{n}_{n-k}(t^{\perp}_{\Omega})$; we have also
$n_{1}(J)=k-1$.
The bijection
$(F^{n}_{k\,n-k}(\Omega))^{-1}$
maps $J$ onto the
the maximal irregular set
satisfying the required conditions.
$\blacksquare$

\subsection{Conclusion}

In the chapter we constructed
three maximal
irregular subsets $I_{1},I_{2}, I_{3}$
of $\ggg^{n}_{k}(V)$
(Proposition 3.1.3 , Theorems 3.2.3 and 3.2.4,
respectively) satisfying the following conditions
$$n_{1}(I_{1})=n_{n-1}(I_{1})=n-k\;,$$
$$n_{1}(I_{2})=n-k-1\;,$$
$$n_{n-1}(I_{3})=n-k+1\;.$$
In the case $n\ne 2k$ these sets are
munually non-similar (see Example 3.2.2).
For the case when $n=2k$ it fails,
the sets $I_{2}$ and $I_{3}$ are similar
(in this case the mapping $F^{n}_{n-k\,k}(\Omega)$
exploited to prove Theorem 3.2.4 is a
regular transformation of $\ggg^{n}_{k}(V)$).
However, the sets $I_{2}$ and $I_{3}$
are not similar.

%% file: BIB.TEX
\thispagestyle{empty}

%% file: IRREG.bbl
\begin{thebibliography}{W99}

\label{bib}
\addtocontents{toc}{\contentsline{chapter}{Bibliography}{\pageref{bib}}}

\bibitem[Art]{Artin} E. Artin, Geometric Algebra, Interscience,
New York, 1957.

\bibitem[BH]{BH}
G. Birkhoff, von J. Neumann,
The logic of quantum mechanics,
Ann. of Math., 37 (1938), 823--843.

\bibitem[Chog]{Chogoshvili}
G. Chogoshvili,
On a theorem in the theory of dimensionality.
Compositio Math., 5 (1938), 292--298.


\bibitem[Chow]{Chow}
W. L. Chow,
On the geometry of algebraic homogeneous spaces,
Ann. of Math., 50 (1949), 32--67.

\bibitem[Die]{Die}
J. Dieudonne,
La Geometrie des groupes classiques, Springer -- Verlag,
Berlin -- New York, 1971.


\bibitem[Dr]{Dr}
A. N. Dranishnikov,
On Chogoshvili's Conjecture,
Proc. Amer. Math. Soc. 125 (1997),
2155--2160.

\bibitem[Eng]{Eng}
R. Engelking, Dimension theory, North -- Holland, 1978.


\bibitem[O'M]{O'M}
O. T. O'Meara, Lectures on linear groups,
Providence, Rhode Island, 1974.

\bibitem[P1]{P1}
M. A. Pankov,
        Projections of $k$-dimensional subsets of $\rr^{n}$
        onto $k$-dimensional planes,
        Matematicheskay fizika, analiz, geometriya, 5
        (1998), v. 1/2, 114--124.

\bibitem[P2]{P2}
M. A. Pankov,
        Projections of $k$-dimensional subsets of $\rr^{n}$
        onto $k$-dimensional planes and irregular subsets of
        the Grassmannian manifolds,
        Topology and its Applications, 101 (2000), 121--135.

\bibitem[P3]{P3}
M. A. Pankov, Irregular subsets of the Grassmannian
manifolds and their mappings,
Matematicheskay fizika, analiz, geometriya, 7 (2000), v. 3.

\bibitem[S]{S}
K. Sitnikov,
An example of a two-dimensional set in
             three-dimensional  Euclidean space allowing
             arbitrarily small deformations into a
             one-dimensional polihedron and a
             certain new characterization of
             the dimension of sets in
             Euclidean spaces,
             Dokl. Akad. Nauk SSSR, 88 (1953), 21--24.

\end{thebibliography}
